\documentclass[a4paper,12pt,tocvesec2]{amsart}

\usepackage{amsmath, amsfonts, amsthm, amsfonts}
\usepackage[francais]{babel}
\usepackage[latin1]{inputenc}

\newcommand{\OOO}{\mathcal{O}}
\newcommand{\X}{\mathcal{X}}
\newcommand{\Ac}{\mathcal{A}}
\newcommand{\Sc}{\mathcal{S}} 
\newcommand{\rhob}{{\bar\rho}}   
\newcommand{\rhobar}{{\bar\rho}}  
\newcommand{\A}{\mathbb{A}}
\newcommand{\ZZ}{\mathbb{Z}}

\newcommand{\Z}{\mathbb{Z}}
\newcommand{\CC}{\mathbb{C}}
\newcommand{\C}{\mathbb{C}}
\newcommand{\RR}{\mathbb{R}}

\newcommand{\R}{\mathbb{R}}
\newcommand{\anneau}{\mathcal{O}}
\newcommand{\QQ}{\mathbb{Q}}
\newcommand{\Q}{\mathbb{Q}}
\newcommand{\Qpb}{\overline{\QQ}_p}

\newcommand{\U}{\textrm{U}}
\newcommand{\NN}{\mathbb{N}}
\newcommand{\Fil}{\textrm{Fil}}
\newcommand{\HH}{\mathcal{H}}
\newcommand{\Gal}{\textrm{\rm{Gal}($\overline{E}/E$)}}
\newcommand{\Dp}{\textrm{Gal($\overline{\QQ}_p/\QQ_p$)}}
\newcommand{\diag}{\textrm{diag}}
\newcommand{\W}{W}

\newcommand{\PGL}{\textrm{PGL}}

\newcommand{\cond}{\textrm{cond}}
\newcommand{\disc}{\textrm{disc}}
\newcommand{\Ker}{\textrm{Ker}}
\newcommand{\Frob}{\textrm{Frob}}
\newcommand{\Hom}{\textrm{Hom}}
\newcommand{\DD}{D_{cris,v_1}}
\newcommand{\Ind}{\textrm{Ind}}
\newcommand{\Spec}{\textrm{Spec}}
\newcommand{\Specmax}{\textrm{Specmax}}
\newcommand{\End}{\textrm{End}}
\newcommand{\tr}{\textrm{tr}}
\newcommand{\Gl}{\textrm{GL}}

\newcommand{\Sl}{\textrm{SL}}
\newcommand{\GL}{\textrm{GL}}
\newcommand{\p}{m}
\newcommand{\sss}{\mathfrak{s}}
\newcommand{\BB}{\mathcal{B}}

\def\]{\textup{\mbox{]\hspace{-.15em}]}}}
\def\[{\textup{\mbox{[\hspace{-.15em}[}}}

\def\got{\mathfrak}

\newtheorem{thm}{Théorème}[section]
\newtheorem{prop}{Proposition}[section]
\newtheorem{cor}{Corollaire}[section]
\newtheorem{lemme}{Lemme}[section]
\newtheorem{remarque}{Remarque}[section]
\newtheorem{conjecture}{Conjecture}[section]

\begin{document}

\title{Formes non temp\'er\'ees pour $\U(3)$ et conjectures de Bloch-Kato}
\author{Joël Bellaïche et Gaëtan Chenevier}
\maketitle

\bigskip

\moveright 0.25in\vbox{\hsize=5in \small
\noindent {\bf Résumé : }
Dans cet article, nous utilisons des familles $p$-adiques de formes
automorphes pour un groupe unitaire à trois variables, passant par des formes 
non tempérées construites par Rogawski, pour montrer certains cas des 
conjectures de Bloch et Kato.

\noindent {\bf Abstract : }
In this paper, we use $p$-adic families of automorphic forms 
for an unitary group in three variables, containing some non-tempered forms constructed by Rogawski, 
to prove some cases of the Bloch-Kato conjectures.
\bigskip}
\bigskip
\tableofcontents

\section{Introduction}
\subsection{Énoncé du théorème principal}

Soit $E/F$ une extension $CM$, $(\rho,V)$ une repré\-sentation de $\Gal$ sur 
une extension finie $L$ de $\Q_p$, continue, irréductible, et géométrique 
(cf. \cite[page 650]{FP}). 
On note $\tau$ l'élément non trivial de $\text{Gal}(E/F)$ et $\rho^\bot$ 
la représentation sur $V^\ast$ donnée par $g \mapsto {}^t\rho(\tau g \tau)^{-1}$.
On suppose que la représentation $\rho$ vérifie\footnote{Notons par exemple
que si $\rho_F$ est une représentation de $\text{Gal}(\bar F/F)$ vérifiant 
$\rho_F^\ast \simeq\rho_F(-1)$, sa restriction $\rho$ à $\Gal$ vérifie~\ref{antiautoduale}. C'est la cas par exemple d'un twist convenable des 
représentation attachés aux formes modulaires}  
\begin{eqnarray} \label{antiautoduale} \rho^\bot = \rho(-1) ,\end{eqnarray}
 où $\rho(-1)$ est un twist à la Tate de $\rho$.

Il est conjecturé que la fonction $L$ complète de $\rho$, notée $L(\rho,s)$,
admet un prolongement méromorphe à tout le plan complexe, holomorphe en zéro si $\rho$ n'est pas le caractère cyclotomique, ce que nous supposerons.
Sous l'hypothèse~\ref{antiautoduale}, cette fonction admet une équation
 fonctionnelle de la forme 
$$L(\rho,s)=\epsilon(\rho,s) L(\rho,-s),\text{ avec }
\epsilon(\rho,0)=\pm 1,$$
En particulier, $\epsilon(\rho,0)=-1$ 
si et seulement si la fonction $L(\rho,s)$ s'annule en $0$ à l'ordre impair.

À la suite de Janssen, Bloch-Kato, Fontaine-Perrin-Riou, notons $H^1_f(\Gal,\rho)$ le sous groupe de $H^1(\Gal,\rho)$ paramétrant les extensions $U$ de la représentation triviale $1$ par $\rho$
$$0 \rightarrow V \rightarrow U \rightarrow 1 \rightarrow 0$$
qui ont bonne réduction en toute place, ce qui signifie que les suites
  $$0 \rightarrow V^{I_w} \rightarrow U^{I_w} \rightarrow 1 \rightarrow 0$$
sont exactes pour toute place $w$ de $E$ ne divisant pas $p$, $I_w$ désignant 
un sous-groupe d'inertie en $w$ et que les suites
$$ 0 \rightarrow D_{cris,w}(V) \rightarrow D_{cris,w}(U) \rightarrow 1 \rightarrow 0$$
sont exactes pour toute place $w$ divisant $p$. La conjecture de Bloch et 
Kato  (cf. \cite[3.4.5]{FP}) implique  que
\begin{conjecture} \label{conjecture}
Si $\epsilon(\rho,0)=-1$, alors $\dim H^1_f(\Gal,\rho) \geq 1$.
\end{conjecture}

Cet article propose une nouvelle méthode pour aborder cette conjecture, basée 
sur des congruences entre formes automorphes non tempérées et tempérées.
Son but est d'en démontrer le cas particulier suivant :

\begin{thm} \label{principal}
Supposons que $E$ est un corps quadratique imaginaire.
Soit $\chi$ un caractère de Hecke algébrique sur $E$, qui vérifie 
$$\chi^\bot=\chi(-1)$$ et dont le type à 
l'infini est de la 
forme
$$z \mapsto z^a \bar z^{1-a},$$
avec $a \geq 2$. Soit $p$ un nombre premier décomposé dans $E$ et non
ramifié pour $\chi$, et $\chi_p:\Gal \rightarrow L^\ast$ une réalisation $p$-adique de $\chi$ sur un corps $L$.
Alors, si $\epsilon(\chi,0)=\epsilon(\chi_p,0)=-1$, on a 
$$\dim H^1_f(\Gal,\chi_p) \geq 1.$$
Autrement dit, il existe une extension non triviale ayant bonne réduction 
partout de la forme 
\begin{eqnarray} \label{extension}
0 \rightarrow \chi_p \rightarrow U \rightarrow 1 \rightarrow 0.
\end{eqnarray}
\end{thm}

Notons que ce théorème peut se démontrer facilement à partir de la conjecture
d'Iwasawa pour les corps quadratiques imaginaires, prouvée par Rubin 
(\cite{rubin}). Cependant, notre méthode est différente (Rubin 
utilise des systèmes 
d'Euler et non des formes automorphes), susceptible de généralisations 
ultérieures (voir plus bas) et donne une information supplémentaire :
elle prouve que {\it les 
réductions modulo $\p^n$, ($\p$ idéal maximal de $\anneau_L$, $n$ entier 
arbitraire) de l'extension~($\ref{extension}$) dont le théorème affirme 
l'existence apparaissent comme sous-quotient de 
la cohomologie étale de variétés algébriques 
sur $E$}. Cet énoncé est prédit par la conjecture de Fontaine-Mazur 
et ne découle pas de la preuve de Rubin.

\subsection{La méthode}

\subsubsection{}
Le théorème~\ref{principal} est une généralisation du théorème principal de 
la thèse d'un des auteurs (cf. \cite{Joelthese} théorème VIII.1.7.2), où il
est prouvé, sous les mêmes hypothèses, l'existence d'une extension non triviale
ayant bonne réduction partout sur $\anneau_L/\p$ de la forme $0 \rightarrow \bar \chi_p \rightarrow U \rightarrow 1 \rightarrow 0$,
et ce pour un ensemble de densité non nulle de $p$.
La méthode utilisée ici est similaire à celle de cette 
thèse, à une variante près (l'emploi de familles $p$-adiques combinées avec un
résultat de Kisin, au lieu d'augmentation du niveau) dont l'idée est due à Urban et Skinner.
Ils ont en effet récemment annoncé (cf. \cite{cras}) un analogue du théorème
\ref{principal} pour des formes modulaires ordinaires de niveau $1$, par une méthode semblable.
Il nous a semblé bon de reprendre cette idée dans notre cas, notamment parce
qu'elle nous permet de se passer d'hypothèses sur $p$, que les familles $p$-adiques
pour $U(n)$ sont construites dans \cite{Che}, et qu'elle est plus simplement généralisable
(voir \ref{generalisations}).
 
\subsubsection{}
Expliquons le principe de la méthode employée. La première idée est d'utiliser
les formules de multiplicités, données par des conjectures d'Arthur, des représentations 
automorphes dans le spectre discret. Pour certaines représentations non 
tempérées, ces formules font apparaître le signe de certains facteurs $\epsilon$. 
Plus précisément, dans le cas du groupe unitaire $\U(3)$ compact à l'infini attaché au corps 
quadratique imaginaire $E$, ces formules sont démontrées par 
Rogawski et affirment que pour $\chi$ un caractère de Hecke comme dans l'énoncé 
du théorème, il existe dans le spectre discret de $E$ une représentation 
$\pi(\chi)$, minimalement ramifiée\footnote{Voir la proposition~\ref{endo-1} 
pour plus de détails}, dont la représentation galoisienne associée 
$\rho:\Gal \rightarrow \Gl_3(L)$ est de la forme 
$\rho=\chi_p \oplus 1 \oplus \chi_p(-1)$
si, et seulement si, $\epsilon(\chi,0)=-1$.

\subsubsection{}
Plaçons-nous sous cette hypothèse. On dispose donc de $\pi(\chi)$ et de sa 
représenta\-tion galoisienne associée $\rho$. L'étape suivante est 
d'obtenir une déformation génériquement irréductible
$\rho'$ de $\rho$, i.e. une représentation $\rho'$ de $\Gal$ sur un 
anneau de valuation discrète $R$ de corps résiduel $L$ 
avec $(\rho'\otimes_R L)^{ss} = \rho$, et $\rho'\otimes_R Frac(R)$ 
irréductible ; il faut également contrôler la ramification de $\rho'$ ainsi 
que son comportement aux places divisant $p$. Le méthode utilisée pour 
construire $\rho'$ consiste à placer $\Pi(\chi)$ dans une famille $p$-adique
de formes automorphes\footnote{C'est là la principale différence 
avec~\cite{Joelthese} où l'on ne construisait qu'une déformation de 
$\rhobar \simeq \bar \chi_p \oplus 1 \oplus \bar \chi_p(1)$ à l'aide d'un 
théorème d'augmentation du niveau. Notons cependant que l'on pourrait
 aussi prouver le théorème~\ref{principal} à l'aide de cette méthode, en 
construisant des déformations de $\rhobar$ qui sont congrues à $\rho$ modulo 
$\p^n$, à l'aide de théorèmes d'augmentation du niveau modulo $\p^n$. On
obtiendrait ainsi des extensions modulo $\p^n$, puis on passerait à la limite sur $n$.
Cette méthode fera l'objet d'un travail ultérieure, dans un cadre 
un peu différent.} 

Pour contrôler la ramification de $\rho'$ aux places $w$ de $E$ 
ne divisant pas $p$, en particulier aux places où $\chi$ est ramifié,
on est obligé d'imposer aux formes de la famille $p$-adique construite de contenir certains
types de Buschnell et Kutsko. Pour traduire l'existence de ces types en termes de la 
ramification de $\rho'$, on bute sur la difficulté suivante : il ne semble
pas connu \footnote{On ne peut appliquer à nos formes le résultat principal 
de \cite{Ha} car elles n'en vérifient pas l'hypothèse, à savoir d'être
de carré intégrable en au moins une place finie.} que la 
construction de la représentation galoisienne $\rho$ attachée à une forme automorphe pour $\U(3)$ est 
compatible, en chaque place, avec la correspondance de Langlands locale.
Nous montrons comment contourner cette difficulté. 
Pour contrôler le comportement de $\rho'$ aux places $w$ divisant $p$,
on utilise une forme convenable d'un résultat récent de Kisin.

\subsubsection{} La dernier étape consiste à appliquer une vieille idée de 
Ribet (cf. \cite{ribet}) : l'existence d'une déformation de $\rho$ comme 
ci-dessus implique l'existence d'extensions non triviales entre les facteurs 
de $\rho$. Dans le cas de Ribet, $\rho$ n'avait que deux facteurs et Ribet 
montrait qu'on pouvait obtenir les extensions d'un facteur par 
l'autre dans les deux sens possibles. Mais c'est un fait incontournable 
(cf. ~\cite{ribellaicheg} )
que quand $\rho$ a plus de deux facteurs irréductibles, on ne peut 
assurer l'existence de toutes les extensions entre ces facteurs. On obtient seulement une disjonction d'assertions d'existence. Autrement dit,
pour montrer  l'existence de l'extension cherchée, on est ramené à montrer la 
non-existence de certaines autres extensions, non-existence qui a une 
signification arithmétique globale, étant aussi un cas
des conjectures de Bloch-Kato.  
Dans~\cite{Joelthese} ainsi que dans la méthode décrite dans~\cite{cras} (voir {\it loc. cit.} dernier paragraphe),
la preuve de ces cas de non-existence repose sur des résultats récents et 
difficiles de Rubin (\cite{rubin}) et de Kato ($p$-adic Hodge theory and values of 
Zeta functions of modular forms, prépublication) respectivement.
Dans cet article, le seul cas de non-existence que nous avons à vérifier
est celui d'une extension ayant bonne réduction partout du 
caractère trivial par le caractère cyclotomique, qui est une application
élémentaire de la théorie de Kummer.

\subsection{Généralisations}\label{generalisations}
Tout d'abord, l'hypothèse sur $p$ dans le théorème \ref{principal} est
inessentielle, et devrait être supprimée lorsque les familles $p$-adiques
seront disponibles aux places non décompos\'ees et en niveau sauvage
(travail en cours d'un des auteurs avec K.Buzzard). \par

 La méthode que nous avons esquissée ci-dessus se prête à des généralisations
aux groupes unitaires $\U(n)$ compact à l'infini associé à un corps 
quadratique imaginaire $E$. Cependant, on ne dispose pas pour l'instant pour ces 
groupes, si $n \geq 4$, de la construction d'une représentation galoisienne 
attachée à chaque forme automorphe et vérifiant les propriétés attendues 
en presque toute place, pas plus que l'on ne dispose des cas nécessaires des 
formules de multiplicités d'Arthur  (cependant des progrès ont été fait récemment dans le cas $n=4$).
Nous avons néanmoins pris garde à énoncer et à démontrer la plupart des lemmes
que nous utilisons sous une forme générale, en vue de leur utilisation 
pour le cas de $\U(n)$. Nous comptons revenir sur ce cas dans un avenir proche.
\par
 Enfin, il semble par contre plus délicat de généraliser cette méthode au 
cas où $E$ est un corps $CM$ quelconque, et non quadratique imaginaire (voir à ce propos la remarque~\ref{corps
CM}). \par 

\subsection{Plan de l'article} 
Indiquons brièvement le contenu des différentes parties de cet article.
Dans la partie 2, nous fixons les conventions (concernant essentiellement 
les normalisations  de la théorie du corps de classes et de la correspondance de Langlands 
locales) et les principales notations que nous utiliserons. 
La partie 3 décrit la construction par Blasius et Rogawski du système de 
représentations $l$-adiques attaché aux formes automorphes pour $U(3)$. 
Nous démontrons quelques résultats concernant la compatibilité de cette construction
avec la correspondance de Langlands locale. La partie 4 construit et décrit
la représentation non tempérée $\pi(\chi)$ discutée plus haut. \par
Les parties 5, 6 et 7 sont rédigées dans une plus grande 
généralité. La partie $5$ fait, pour $\textrm{GL(n)}$, la
théorie des différents choix possibles d'une forme propre ancienne de niveau 
iwahorique en $p$, attachée à une forme non ramifiée en $p$. La partie $6$ 
traite des raffinements des représentations $p$-adiques cristallines, 
contrepartie galoisienne de la théorie précédente, et des 
familles de telles représentations raffinées. On y énonce en particulier,
sous une forme adaptée à notre usage, un résultat récent de Kisin. 
La partie 7 contient les résultats nécessaires sur les pseudo-caractères et les 
représentations galoisiennes, ainsi que la généralisation adéquate du lemme 
de Ribet.

Dans la partie 8 nous revenons aux groupes unitaires à trois variables,
et construisons des familles $p$-adiques passant par $\pi(\chi)$. Enfin, la partie 9 montre le 
théorème principal. \par \vspace{1 mm}
 
{\small Les auteurs\footnote{Durant la rédaction de cet article l'un des 
auteurs (Joël Bellaïche) a bénéficié de l'aide du Bell fund et du Ellentuck 
Fund} sont heureux de remercier Laurent Berger, Pierre Colmez,
Laurent Clozel, Michael Harris, Guy Henniart et Éric Urban
pour de nombreuses et utiles conversations durant la réalisation de cet 
article.}

\section{Notations et conventions}\label{notations}

\subsection{Corps et groupes de Galois} Dans tout cet article, on note $E \subset \CC$ une corps quadratique 
imaginaire, $\bar E$ la clôture algébrique de $E$ dans $\CC$, 
$\A_E$ l'anneau des adèles de $E$, $\W_E$ et $\Gal$
les groupes de Weil et de Galois de $\bar E$ sur $E$. Pour $v$ place
de $E$ on notera $D_v$ un sous-groupe de décomposition de $\Gal$ en
$v$, $I_v$ le sous-groupe d'inertie de $D_v$, et $\Frob_v \in D_v$ un
Frobénius géométrique. On note $\tau \in
\text{Gal}(\bar E/\QQ)$ la conjugaison complexe, $\tau^2=1$. 

\subsection{Anti-auto-duale} Soit $\rho$ une représentation de $\W_E$ ou de $\Gal$ dans
$\GL_n(A)$, où $A$ est un anneau commutatif. On notera $\rho^\bot$ la
représentation définie par 
$$\rho^\bot(g)={}^t \rho(\tau g \tau)^{-1}$$ 
\par \vspace{1 mm}

\subsection{Représentations algébriques} \label{plushautpoids}Soit $w=(k_1 \geq \cdots \geq k_n) \in \ZZ^n$, on note $V_w$ la représentation algébrique de
$\GL_n/\QQ$ de plus haut poids $w$ relativement au Borel supérieur de $\GL_n$. C'est l'unique représentation algébrique irréductible de
$\GL_n/\QQ$ telle que le tore diagonal de $\GL_n$ 
agissent sur la droite stable par le Borel supérieur
par $\diag(x_1,\cdots,x_n) \mapsto \prod_{i=1}^n x_i^{k_i}$. On note
$$\ZZ^{n,+}:=\{(k_1,...,k_n) \in \ZZ^n, k_1\geq ... \geq k_n \}$$ 
$$\ZZ^{n,--}:=\{(k_1,...,k_n) \in \ZZ^n, k_1<k_2<...<k_n\}$$
Si $w=(k_1,...,k_n) \in \ZZ^{n}$, on pose $$-w=(-k_n,...,-k_1) \in
\ZZ^{n} \, \, \textrm{ et  } \, \, \, \delta(w)=\textrm{Min}_{i=1}^{n-1} (k_i-k_{i+1})
\in \ZZ$$  Si $w \in \ZZ^{n,+}$, $-w \in \ZZ^{n,+}$ et $\delta(w) \in \NN$. 
Le dual de $V_w$ est alors $V_{-w}$. Si $F$ est un
corps de caractéristique $0$, on notera $V_w(F)$ la représentation de
$\GL_n/F$ extension des scalaires à $F$ de $V_w$. \par \vspace{2 mm}

\subsection{Correspondance de Langlands locale} Nous précisons dans ce paragraphe les conventions choisies
pour fixer l'isomorphisme de réciprocité de la théorie du corps de classes ainsi que la correspondance
de Langlands locale pour $\GL_n$. Soit $F$ un corps de nombres, on normalise l'isomorphisme d'Artin de la
théorie du corps de classes globale
$$\textrm{rec}_F :\A_F^\ast / \overline{F^\ast(F\otimes_{\QQ}\RR)^{0,*}} \rightarrow
\textrm{Gal}(\bar F/F)^{ab}$$
 en demandant qu'il envoie toute uniformisante locale $\pi_v$ en une place
finie $v$ sur le Frobenius {\it géométrique} de $D_v/I_v$, 
avec l'abus de langage évident. Pour toute place $v$ de $F$, on dispose
alors par restriction de $\textrm{rec}_F$, d'un isomorphisme du corps de
classes local
$$\textrm{rec}_{F_v}\, \, : \, \,  F_v^\ast \longrightarrow \W_F^{ab}$$ compatible à l'isomorphisme global.
\par \vspace{1 mm}
Supposons maintenant que $F$ est un corps local non archimédien. On choisit
la normalisation {\it à la Langlands} pour la correspondance de Langlands
locale pour $\GL_n(F)$, notée $L$ (voir par exemple \cite{Ha} p.2). Ainsi, pour $\pi$ une
représentation irréductible lisse de $\GL_n(F)$, $L(\pi)$ est une
représentation complexe $\Phi$-semisimple de dimension $n$ du groupe de Weil-Deligne
$\W \!D_F$ (cf. \cite[4.1]{Tate}).
Pour $n=1$, $\pi$ est un caractère de $GL_1(F)=F^\ast$ et $L(\pi)$
est le caractère de $\W_F$ qui s'en déduit via l'isomorphisme
$\textrm{rec}_{F_v}$ ci-dessus. Par exemple, $L$ envoie le quotient de Langlands de l'induite 
parabolique attachée à deux représentations lisses irréductible
$\pi_1$ et $\pi_2$ sur la somme directe $L(\pi_1) \oplus L(\pi_2)$.

Soit $l$ un nombre premier, on fixe $\iota_l: \CC \rightarrow \overline{\QQ}_l$
un isomorphisme de corps. Pour $\pi$ une représentation lisse irréductible
de $\GL_n(F)$, $F$ un corps local non archimédien, on peut voir $\L(\pi)$, par transport 
de structure via $\iota_l$ comme une représentation sur $\overline{\QQ}_l$ en fait 
définie sur une extension finie de $\Q_l$, et on peut donc lui associer une 
représentation notée $L_l(\pi)$ du groupe de Weil ordinaire $\W_F$
sur $\bar \Q_l$, comme en \cite[4.2.1]{Tate}.

\subsection{Caractères de Hecke} Soit $\chi: E^*\backslash \A_E^* \rightarrow \CC^*$ 
un caractère de Hecke de $E$, pour toute place $v$ de $E$ on note $\chi_v$
la restriction de $\chi$ à $E_v^*$ et $\chi_f$ sa restriction aux idèles
finies $\A_{f,E}^*$. On identifie $E\otimes_{\QQ} \RR$ à $\CC$ via
l'inclusion $E \subset \CC$. On suppose dans ce qui suit que $\chi$ est algébrique, i.e
$\chi_{\infty}(z)=z^a \bar{z}^b$ avec $a,b \in \ZZ$, et on fixe encore
$l$ et $\iota_l$ comme dans la section précédente en supposant de plus
$l=v_1v_2$ totalement décomposé dans $E$, $v_1$ étant la place définie par $E
\subset \CC \overset{\iota_l}{\rightarrow} \overline{\QQ_l}$. À
$(\chi,\iota_l)$ on
peut alors associer un caractère continu $\overline{E^*}\backslash \A_{f,E}^* \rightarrow
\overline{\QQ}_l^*$ défini par la formule: 
$$x \mapsto \iota_l(\chi_f(x)))x_{v_1}^ax_{v_2}^b$$
Par composition avec $\textrm{rec}_E$, on en déduit un caractère
$l$-adique de $\Gal$ que nous noterons $\chi_l$ (à ne pas
confondre avec une composante locale de $\chi$, mais c'est sans ambiguïté). 

\subsection{Poids de Hodge} \label{notationshodge} 
Dans ce paragraphe, $p$ est un nombre 
premier fixé,
$\mu_{p^n}$ désigne le $\Dp$-module des racines $p^n$-ièmes de l'unité de
$\overline{\QQ}_p$. On note $\omega$ le caractère cyclotomique $\Dp
\rightarrow \End_{\ZZ_p}(\underset{ n \leftarrow \infty}{\lim}
\mu_{p^n})=\ZZ_p^*$. On note $\QQ_p(1)$ le $\QQ_p$-espace vectoriel de dimension $1$ muni d'une action de $\Dp$ par le
caractère cyclotomique. Notre convention est que cette représentation est 
de Hodge-Tate de poids de Hodge $-1$ (cf. \cite{Sen1}). Si $n \in \ZZ$, on note
$\QQ_p(n):=\QQ_p(1)^{\otimes n}$, et si $V$ est un $\QQ_p$-espace vectoriel 
représentation
de $\Dp$, $V(n):=V \otimes_{\QQ_p} \QQ_p(n)$. \par \vspace{1 mm} 
Si $F$ est un corps local et $V$ un $F$-espace vectoriel de dimension finie 
qui est une représentation continue de $\Dp$, nous noterons
$D_{cris}(V):=(V\otimes_{\QQ_p} B_{cris})^{\Dp}$ (cf. \cite{Fo1} \S 2.3,
3.1). Il hérite de $B_{cris}$ d'un endomorphisme $K$-linéaire $\varphi$, le Frobénius cristallin, 
et lorsque nous parlerons des valeurs propres de ce dernier, ce sera toujours vu comme $K$-endomorphisme.
On rappelle que $V$ est dite cristalline sur $\dim_{K}D_{cris}(V)=\dim_{K}(V)$
(cf. \cite{Fo2} 5.4). On parlera parfois, par abus, du Frobénius
cristallin de $V$ pour celui de $D_{cris}(V)$.\par \vspace{1 mm}

\subsection{Géométrie rigide} \label{notationsrigides} Si $X/F$ est un affinoide sur un
corps local $F$, on notera $A(X)$ la $F$-algèbre affinoide de $X$. Si $X$
est réduit, la norme sup. sur $X$ fait de $A(X)$ une $F$-algèbre de Banach
commutative. On notera $A(X)^0$ les éléments de $A(X)$ de norme $\leq 1$.
$\A^n$ désigne l'espace affine rigide analytique de dimension $n$ sur
$\QQ_p$.

\section{Rappel sur la classification de Rogawski}

\subsection{Les groupes unitaires considérés} Soit $f$ la forme hermitienne sur
$E^3$ de matrice 
$$ \left( \begin{array}{ccc} & & -1 \\ & 1 &  \\ -1 & & \end{array} \right )$$

Comme Rogawski (\cite{Rog1} p. 66,67), on note $\U(2,1)$ le groupe
unitaire sur $\ZZ$ défini par cette forme, il est quasi-déployé.
On fixe $\U(3)$ une forme intérieure de $\U(2,1)$ compacte à l'infini, 
inchangée aux places finies (\cite[lemme 2.1]{Clo}). 

\subsection{Classification de Rogawski} \label{classification}

\subsubsection{}
Suivant Rogawski, et comme prédit par les conjectures d'Arthur, 
les représentations automorphes discrètes des groupes $\U(2,1)$ et
$\U(3)$ sont regroupées en $A$-paquets de $5$ types (\cite[2.9]{Rog1}). 
Chaque $A$-paquet $\Pi$ possède un changement de base $\pi_E$ qui est une représentation
automorphe de $\GL_3/E$ (\cite[2.8]{Rog1}). 

\subsubsection{} \label{propgal} Étant donné un nombre premier $l$ et un isomorphisme de corps
$\iota_l: \CC \rightarrow \overline{\QQ}_l$, on peut associer grâce au travaux de Rogawski
à un $A$-paquet $\Pi$ de $\U(3)$  (resp. de $\U(2,1)$ s'il est cohomologique à l'infini) une représentation
$l$-adique continue, semi-simple, découpée dans la cohomologie $l$-adique
des surfaces de Picard:
$$\rho_l(\Pi) : \Gal \rightarrow \GL_3(\overline{\QQ}_{l}),$$
caractérisée par la propriété suivante, pour toute place finie $w$ de $E$: 
\begin{eqnarray} \label{etoile}
\text{   Si $w$ ne divise pas $l\disc(E)$ et si $(\pi_E)_w$ est non ramifiée,  }
\rho_l(\Pi)_{|D_w} \simeq  L_l((\pi_E)_w)
\end{eqnarray}

Autrement dit $\rho_l$ est non ramifiée en $w$, et le polynôme
caractéristique d'un Frobénius géométrique $\Frob_w$ est égal à celui de la
matrice de Hecke de $(\pi_E)_w$. En conséquence, $\rho_l \simeq
\rho_l^\bot$. Lorsque $\Pi$ est sous-entendu, on notera $\rho_l$ pour
$\rho_l(\Pi)$. Bien que le choix de $\iota_l$ n'est pas apparent dans la
notation $\rho_l$, il sera toujours sous-entendu.

\subsubsection{} Nous nommons et décrivons ci-dessous les 5 types de $A$-paquets, ainsi que 
les propriétés des représentations galoisiennes associées quand elles
existent. Quand le $A$-paramètre $a$ et
le $L$-paramètre $\phi$ du changement de base $\pi_E$ du paquet considéré ont un sens non 
conjectural (i.e. se factorise par
le quotient $W_E \times SL_2(\CC)$ du groupe conjectural $L_E \times
SL_2(\CC)$), nous les donnons. \par \vspace{2 mm}
	L'existence de $\rho_l$ satisfaisant \ref{etoile} résulte, dans les
cas non triviaux, du théorème 1.9.1 de \cite{br1}, mais le
lecteur prendra garde que la définition de $\rho_l$ que nous
prenons (motivée par la  vérification de~\ref{etoile} \S \ref{propgal}) ne coïncide pas avec la
représentation appelée $\rho_l$ {\it loc. cit.}, que nous noterons 
$\rho_{l,rog}$ ci-dessous (qui d'ailleurs n'est pas toujours de dimension 3).

\vspace{2 mm} \begin{itemize} \par \vspace{1 mm}
\item Cas stable tempéré; l'existence de $\rho_l$
résulte de~\cite[théorème 1.9.1 (a)]{br1}. Notre $\rho_l$ est 
$\rho_{l,rog}(1)$ avec les notations de {\it loc. cit}. La repr\'esentation 
$\rho_l$ est irréductible\footnote{Nous ne nous servirons pas de ce fait. Voir la remarque suivant la proposition~\ref{irrgen}} et satisfait $\rho_l
\simeq \rho_l^{\bot}$. \par \vspace{2 mm}

\item Cas endoscopique tempéré de type $(2,1)$;  l'existence de 
$\rho_l$ résulte de \cite[théorème 1.9.1 (b)]{br1} (on définit $\rho_l$
comme $\rho_{l,rog}(1) \otimes \chi_l^{-1} \oplus \chi_l(1)$
avec les notations de {\it loc. cit.}).

On a $\rho_l \simeq \tau_l \oplus \chi_l$, 
avec $\tau_l$ irréductible de dimension $2$, 
$\tau_l \simeq \tau_l^{\bot}$, $\chi_l=\chi_l^{\bot}$. \par \vspace{2 mm}

\item Cas endoscopique tempéré de type $(1,1,1)$ ; le $A$-paramètre 
$a$ est trivial sur facteur $SL_2(\CC)$, et l'on a  $$a_{|L_E}=\phi=\psi_1
\oplus \psi_2 \oplus \psi_3,$$ $\psi_i : W_E \rightarrow \CC^*$ vérifiant
$\psi_i=\psi_i^{\bot}$. Les $A$-paquets de ce type sont cohomologiques 
quand les caractères de Hecke $\psi_i$ sont algébriques, on définit $\rho_l$ par
$(\psi_1)_l \oplus (\psi_2)_l \oplus (\psi_3)_l$ \par \vspace{2 mm}
 
\item Cas endoscopique non tempéré; le $A$-paramètre $a$ vérifie
$$\forall w \in \W_E,\ a(w)=\left( \begin{array}{ccc} \chi(w) & & \\ &
\psi(w) & \\ & &
\chi(w) \end{array} \right)$$
où $\chi$ et $\psi$ sont des caractères de Hecke de $E$ vérifiant $\chi=\chi^\bot$, $\psi=\psi^\bot$,
et $$a_{|\Sl_2(\CC)}\left( \begin{array}{cc} \alpha & \beta \\ \gamma & \delta 
\end{array}\right) = \left(\begin{array}{ccc} \alpha & 0 & \beta \\ 0
& 1 & 0 \\ \gamma & 0 & \delta \end{array} \right).$$
Nous noterons le $A$-paquets correspondant $\Pi(\chi,\psi)$.
Le $L$-paramètre $\phi$ vérifie donc 
$$\forall w \in \W_E,\ \phi(w)=\left( \begin{array}{ccc} \chi(w) |w|^{1/2} & 
& \\ & \psi(w) & \\ & &
\chi(w)|w|^{-1/2} \end{array} \right),$$ et l'on pose, si $\chi
|\ |^{1/2}$ et $\psi$ sont algébriques, auquel cas $\phi$ est
cohomologique à l'infini $$\rho_l= (\chi |\  |^{1/2})_l \oplus \psi_l \oplus  (\chi |\ |^{-1/2})_l.$$
Nous analyserons les $A$-paquets de ce type de manière beaucoup plus
détaillée dans la section suivante. \par \vspace{2 mm}

\item Cas stable non tempéré ; le $A$-paramètre $a$ vérifie
$$\forall w \in \W_E,\ a(w)=\left( \begin{array}{ccc} \chi(w) & & \\ &
\chi(w) & \\ & &
\chi(w) \end{array} \right),$$
où $\chi$ est un caractère de Hecke de $E$ vérifiant $\chi=\chi^\bot$
et $a_{|\Sl_2(\CC)}$ est la représentation irréductible de dimension
$3$.
On a donc pour $L$-paramètre $$\phi(w)=\left( \begin{array}{ccc}
\chi(w)|.|^{-1} & & \\ & \chi(w) & \\ & &
\chi(w)|.| \end{array} \right),$$
où $\chi \simeq \chi^\bot$.
Les $A$-paquets correspondants sont des singletons, 
composés des représentations automorphes de dimension $1$.
\par \vspace{1 mm}
\end{itemize}
\subsection{Propriétés de $\rho_l$.}

On conjecture naturellement que la construction $\Pi \mapsto \rho_l$
est compatible avec la correspondance de Langlands locale. Plus précisément,
on s'attend à ce que l'énoncé suivant soit vérifié: \par \vspace{2 mm}

(COMP) {\it Si $\Pi$ est un $A$-paquet, de changement de base $\pi_E$,
alors pour tout $l$, et pour toute place finie $v$ de $E$ ne divisant pas
$l$, la Frobenius-semi-simplifié de la représentation galoisienne
$(\rho_l)_{|D_v}$ est isomorphe à $L_l((\pi_E)_v)$.} \par \vspace{2 mm}
Il est vrai pour les places finies $v$ ne divisant pas $\disc(E)\cond(\chi_0)$, (c'est la
propriété~\ref{etoile}, \S \ref{propgal}). De plus, 
on pourrait montrer qu'il est vrai pour tout $v$ si $\pi_E$ est de carré
intégrable à au moins une place finie, en utilisant \cite[théorème ]{Ha}. Malheureusement cette hypothèse ne sera jamais vérifiée dans les cas que
nous aurons à traiter.  \par \vspace{1 mm} 
	Nous allons démontrer dans ce qui suit les cas particuliers
de (COMP) dont nous aurons besoin. La propriété suivante est une extension de la propriété~\ref{etoile} aux
places divisant $\disc(E)$. 
Il sera utile de nous placer dans un contexte un peu plus général :
Soit $F$ un corps totalement réel, $\U(3)/F$ le groupe unitaire
à trois variables associé à l'extension $EF/F$, compact à toutes les
places à l'infini. Pour $\Pi$ un $A$-paquet de $\U(3)/F$, les travaux de Rogawski définissent encore, tout comme
dans le cas  $F=\Q$, un changement de base $\pi_{EF}$ à $\Gl_3/EF$, et des représentations galoisiennes 
$\rho_l$ de $\textrm{Gal}(\overline{\QQ}/EF)$, vérifiant l'analogue de ~\ref{etoile} \S \ref{propgal}.  

\begin{prop}\label{langram}
Si $\Pi$ est un $A$-paquet de $\U(3)_F$, de changement de base $\pi_{EF}$,
alors pour tout $l$, et pour toute place finie $v$ de $EF$ ne divisant pas
$l$ où $\pi_{EF}$ est non ramifiée,
$\rho_l{|D_v}$ est 
non ramifiée. \end{prop}

{\it Preuve:} 
Soit $v$ une place finie de $EF$ ne divisant pas $l$ et telle que
$(\pi_{EF})_v$ est non ramifiée, on veut montrer que $\rho_l$ est non ramifiée
en $v$. On note $w$ la place de $F$ que divise $v$. Montrons d'abord le
résultat si $w$ ne divise pas le discriminant relatif de $EF/F$. \par 

Soit $G'_F$ désigne la forme intérieure de $\U(3)_F$ quasi-déployée à 
toutes les places finies de $F$ et à une seule place
infinie, par~\cite[théorème 2.6.1]{Rog1}, le $A$-paquet $\Pi$
se transfère en un $A$-paquet $\Pi'$ cohomologique du groupe $G'_F$. $\rho_l$ 
est alors construite dans la cohomologie de la variété de Shimura
associé à $G'_F$, sans structure de niveau en $w$. Or cette variété de 
Shimura a bonne réduction en $v$ quand $v$ est non ramifié dans $EF/F$ :
cela résulte d'un simple calcul de déformation du problèmes de modules dont
cette variété de Shimura est l'espace de module, comme celui fait 
en~\cite[proposition I.2.1.5]{Joelthese}.
Cela conclut. \par
	On suppose donc que $v$ divise le discriminant de $EF/F$, en
particulier $v$ divise un premier $p \in \ZZ$ ramifié dans $E$, $p \neq l$. Soit $F'$ un corps quadratique réel ramifié en $p$ et tel que
$EF'/F'$ est décomposé en l'unique place de $F'$ au-dessus de $p$. Il est aisé
de voir qu'un tel corps existe toujours. Par exemple si $p\neq 2$,
$E=\Q(\sqrt{pD})$ avec $-D \in \NN$ sans facteur carré, alors il existe $D' \in \NN$
premier à $p$ et sans facteur carré tel que $DD'$ est un carré modulo $p$,
$F'=\QQ(\sqrt{pD'})$ convient. Alors $EF$ ne contient pas $F'$, $FF'$ est totalement réel 
et ramifié au dessus de $w$, en une place que l'on note $w'$, et $EFF'/FF'$ est alors décomposé au dessus de $w'$. \par 
Notons $\pi_{EFF'}$ le changement de base, à la Arthur-Clozel (\cite{AC}, théorème 4.2, chapitre 3) de
$\pi_{EF}$ à ${\Gl_3}_{EFF'}$. La représentation $\pi_{EFF'}$ satisfait
$\pi_{EFF'}^{\tau} \simeq \pi_{EFF'}^*$, car c'est déjà le cas de $\pi_{EF}$. D'après Rogawski (\cite[page xi]{Roglivre})
elle descend donc en un $A$-paquet $\Pi_{FF'}$ du groupe unitaire à trois variables
quasi-déployé $\U(2,1)_{FF'}$ sur $FF'$, et même à $\U(3)_{FF'}$. Mais s'il on suit les
constructions précédentes, un groupe de décomposition en $v$ de $\textrm{Gal}(\bar \Q/EF)$
s'identifie à un groupe de décomposition en $w'$ de $\textrm{Gal}(\bar \Q/EFF')$, et
$\pi_{EFF'}$ est toujours non ramifiée en $w'$. On applique alors le premier
paragraphe à $FF'$ et $\Pi_{FF'}$, ce qui conclut. $\square$ \par \vspace{1 mm}

\begin{prop} \label{pastropram}  Soit $\Pi$ un $A$-paquet de $\U(3)$,
$l$ un nombre premier, $v$ une place finie de $E$ ne divisant pas $l$, telle que
$L((\pi_{E})_v)= \phi_1 \oplus \phi_2 \oplus \phi_3$, où les $\phi_i$ sont des
caractères continus $E_v^* \rightarrow \C^*$. 
Soit $I':=\textrm{rec}_v(\Ker ( (\phi_1 \oplus \phi_2  \oplus
\phi_3)_{|\mathcal{O}_{E_v}^*})) $, alors $\rho_l(I')=1$.
\end{prop}
{\it Preuve:} Soit $p \in \ZZ$ premier au dessous de $v$, considérons 
l'extension abélienne finie de $E_v$ définie par 
$$M:=\overline{E_v}^{\textrm{rec}_{E_v}(\Ker(\phi_1 \oplus \phi_2  \oplus
\phi_3))}$$ 
Supposons que l'on sache trouver un corps de nombres $F/\QQ$ ayant les propriétés suivantes :
\begin{itemize}
\item[i.] $F$ est totalement réel,
\item[ii.] $F$ est construit à partir de $\QQ$ et d'extensions abéliennes successives, 
\item[iii] $EF$ admet une place $w$ divisant $v$ telle que $(EF)_w/E_v$ soit
isomorphe à $M/E_v$.
\end{itemize}
\par 
\vspace{1 mm}

Soit $\pi_{EF}$ le changement de base à $EF$ de $\pi_E$ (\cite[théorème 4.2 
et 5.1]{AC}, applicable par ii.), 
par compatibilité du changement de base global au changement de base local appliqué à l'extension
$(EF)_w/E_v$, on a que $(\pi_{EF})_w$ est non ramifiée (utilisant iii.). 
Si $\Pi_F$ est le $A$-paquet de $\U(3)_F$ changement de base à $F$ de 
$\Pi$ par le changement de base de Rogawski, son changement de base à $E$
est $\pi_{EF}$ par multiplicité $1$ forte dans le spectre discret de $\GL(3)$. La
proposition~\ref{langram} conclut. \par 
	Il reste à trouver un tel corps de nombres $F$. 
Fixons $K/\QQ$ un corps quadratique réel ayant une place $v'$
divisant $p$ tel que $K_{v'} \simeq_{\QQ_p} E_v$. D'après Artin-Tate 
(\cite[théorème 5, page 103]{at}), on
peut trouver une extension abélienne finie $F$ de $K$ totalement réelle, et ayant
une place $v''$ divisant $v'$ telle que $F_{v''}/K_{v'}$ soit isomorphe à
$M/E_v$. $EF \simeq E \otimes_{\QQ} F$ est alors décomposé au dessus de $F_{v''}$, le choix d'une
place $w$ de $EF$ au dessus de $v''$ satisfait donc iii, ce qui conclut. $\square$ \par \vspace{2 mm}

\begin{prop} \label{nonramcris} Soit $l$ un nombre premier décomposé
dans $E$,  $\iota_l$ comme plus haut, $\Pi$ un $A$-paquet pour $\U(3)$ de changement de base $\pi_E$, $v_1$ une
place de $E$ divisant $l$ telle que $(\pi_E)_{v}$ est non ramifiée, 
alors : \begin{enumerate}
\item  $(\rho_l(\Pi))_{|D_{v}}$ est cristalline en $v$.
\par \vspace{1 mm}
\item Le polynôme caractéristique du Frobénius cristallin de $(\rho_l(\Pi))_{|D_{v}}$ est l'image par $\iota_l$
de $L((\pi_E)_{v})(\Frob_{v})$, $\Frob_{v}$ étant un Frobénius géométrique de
$\W_{F_{v}}$,  
\item Si le $L$-paramètre de $\Pi_{\infty}$ a sa restriction à
$\W_{\CC}=\CC^*$ de la forme $(z/\bar{z})^{a_1} \oplus (z/\bar{z})^{a_2}
\oplus (z/\bar{z})^{a_3}$ avec $(a_1,a_2,a_3) \in \ZZ^3$, $a_1> a_2 > a_3$,
les poids de Hodge-Tate de $(\rho_l(\Pi))_{|D_{v}}$ sont les $a_i$ si $v$
est la place de $E$ induite par $E \subset \CC \overset{\iota_l}\rightarrow
\overline{\QQ}_l$, les $-a_i$ sinon.
\end{enumerate}
\end{prop} 

{\it Preuve:} La proposition est immédiate si le $L$-paramètre de $\Pi$ est
somme de trois caractères, il ne reste donc qu'à traiter les cas où $\Pi$ est
stable tempéré ou endoscopique tempéré de type $(2,1)$ (cf. \S \ref{classification}).
Nous nous placerons par exemple dans le premier cas, le second se traitant de
manière identique. \par 
	Soit $F$ un corps quadratique réel décomposé au dessus de $l$, 
$\Pi_F$ le changement de base à $F$ du $A$-paquet $\Pi$ (\cite{Rog1} 2.6.1),
d'après \cite[page xi]{Roglivre} 
$\Pi_F$ descend en un $A$-paquet cohomologique du groupe $H/F$ suivant : $H/F$ est une
forme intérieure du groupe quasi-déployé sur $\U(2,1)/F$, non quasi-déployée
uniquement à l'une des deux places infinies. Supposons donc que $\Pi_F$ est
stable tempéré, on peut choisir une représentation non ramifiée en $l$,
disons $\pi$, dans ce $A$-paquet. 
Par \cite{Rog1} \S 4.3, $\rho_l(\pi)(-1)$ apparaît dans la cohomologie
$l$-adique d'un schéma abélien $A$ sur une certaine surface de Picard $X_U/O_{EF}$,
cette dernière étant projective sur $O_{EF}$ par l'hypothèse sur $H$, 
et lisse sur $O_{(EF)_v}$ car on peut la choisir 
un niveau $U$ qui soit maximal hyperspécial en $v$ par notre hypothèse sur 
$v$. L'assertion $1$ résulte alors de \cite{Cris} . \par 
	Il est bien connu qu'alors la propriété $2$ découle de la construction de
motivique de $\rho_l(\Pi)(-1)$ (\cite{Rog1} \S 4.3) et d'un théorème de Katz-Messing 
(\cite[théorème 2.2]{KMe}), combinés à (\ref{etoile}) \S \ref{propgal} et au changement de base lisse en
cohomologie $l$-adique. De même, la propriété $3$ découle d'un théorème de G.Faltings (\cite{Fa}).
$\square$. \par

\section{Représentation non tempérée attachée à un caractère de Hecke}

\subsection{Hypothèses sur le caractère de Hecke} 
Soit $\chi_0$ un caractère de Hecke de $E$, vérifiant 
$$\chi_0(z\bar{z})=1, \, \forall z \in \A_E^*, \, \, \, \, \, \, \,
\textrm{ et \,  \, } 
\chi_{0,\infty}(z_{\infty})=z_{\infty}^{k}/|z_\infty|^{k}, \textrm{
$k$ est un  entier impair positif } $$

Par la première hypothèse, sa fonction $L$ complète $L(\chi_0,s)$ a pour équation fonctionnelle
(\cite{Tate} 3.6.8 et 3.6.1)
$$L(\chi_0,s)=\varepsilon(\chi_0,s) L(\chi_0,1-s), \, \, \,
\varepsilon(\chi_0,1/2)=\pm 1 $$
On utilisera de plus le caractère algébrique $\chi:=\chi_0 |.|^{1/2}$, dont on note
$\chi_v$ la composante locale en une place $v$ de $E$. On a
$\chi^{\bot}=\chi |.|^{-1}$, et si $\infty$ désigne $E \subset \CC$, 
$\chi_{\infty}(z)=z^{(k+1)/2}\bar{z}^{(1-k)/2}$. On notera
$\cond(\chi_0)$ le conducteur de $\chi_0$, qui est aussi celui de $\chi$, et $\disc(E)$ le discriminant de $E$.
\par \vspace{2 mm}

Dans les sous-parties qui suivent, nous décrivons, suivant Rogawski,
 les composantes locales du $A$-paquet endoscopique non tempéré 
$\Pi(\chi,1)$.

\subsection{Composantes locales en $p$ décomposé} 

\subsubsection{} \label{pdeccan}
Soit $p=v_1v_2$ est un nombre premier décomposé dans $E$, la donnée des
$v_i$ nous fournit un isomorphisme $\U(3)(\QQ_p) \rightarrow
\GL_3(E_{v_i})=\GL_3(\QQ_p)$. Il faut noter que ces deux isomorphismes
diffèrent d'un automorphisme extérieur de $\GL_3(\QQ_p)$. Fixons celui avec
$v_1$ par exemple. \par 

\subsubsection{} \label{descpidec} Soit $P=MN$ le parabolique standard de
$\U(3)(\Q_p)=\GL_3(\Q_p)$ de type $(2,1)$, considérons le 
caractère complexe lisse $\lambda_p$ de $M=\GL_2(\Q_p) \times \GL_1(\Q_p)$ défini par 
$\lambda_p(x,y)=\chi_{0,v_1}(\det(x))$. L'induite parabolique normalisée de 
$\lambda_p$ est irréductible (\cite{Rog2} p.196), on la note
$\pi_p^n(\chi_0)$. \par \vspace{2 mm}

\subsubsection{} \label{kpdecompose}
Si $\chi$ est non ramifié en $p$ (décomposé), $\pi_p^n(\chi_0)$ est non ramifiée, et a donc un vecteur fixe par n'importe
quel compact maximal. On en choisit un, noté $K_p$, égal à $\U(3)(\Z_p)$
pour $p$ assez grand (\cite[3.9.1]{tits}). 

\subsection{Composantes locales en $p$ inerte ou ramifiée}

\subsubsection{} \label{descpiinerte}

Suivant Rogawski (\cite{Rog2} p. 396), définissons pour $p$ un nombre 
premier inerte ou ramifié dans $F$,  une représentation lisse irréductible 
$\pi_p^n(\chi_0)$ de $\U(2,1)(\Q_p)=\U(3)(\Q_p)$ de la manière suivante :

On voit $\chi_{0,p}$ comme caractère du tore diagonal de $\U(3)(\Q_p)$ par :
$$ \left( \begin{array}{ccc} \alpha & &  \\ & \beta &  \\ & &
\bar{\alpha}^{-1} \end{array} \right ) \mapsto \chi_{0,p}(\alpha) $$
L'induite parabolique normalisée de $\chi_p$ a deux facteurs de
Jordan-Hölder, dont l'un est non tempéré, on le note $\pi_p^n(\chi_0)$
(\cite{Roglivre} p. 173, \cite{KE} p. 126) \par \vspace{1 mm}

\subsubsection{}
\label{kpinerte}
Si $p$ est inerte et $\chi$ non ramifié en $p$, $\pi_p^n(\chi_0)$ est non ramifiée et a donc un vecteur fixe par n'importe quel compact maximal 
hyperspécial. Là encore, on en choisit un $K_p$, que l'on prend égal à 
$\U(3)(\Z_p)$ quand on le peut (pour presque tout $p$, cf. \cite[3.9.1]{tits})

\subsubsection{}
\label{kpramifiee} 
Si $p$ est ramifié, le groupe $\U(3)(\Q_p)$ 
est un groupe unitaire ramifié et n'a donc pas de compact maximal 
hyperspécial, mais l'une des deux classes de conjugaison de compacts maximaux
est {\it très spéciale} au sens de \cite[page 88]{labesse}. Si $\chi$ est 
non ramifié en $p$, la représentation $\pi_p^n(\chi_0)$ admet un vecteur fixe par n'importe quel compact maximal très spécial 
d'après~\cite[prop 3.6.2]{labesse}. On en choisit un que l'on note $K_p$.

\subsection{Composantes locales à l'infini}

Tout $L$-paramètre "relevant" de $\U(3)(\R)$ a sa restriction à $W_{\C}=\CC^*$
de la forme :
$$ z\mapsto (z/{\bar z})^{a_1} \oplus {(z/ {\bar z})}^{a_2} \oplus  {(z/{\bar
z})}^{a_3}, \, \, \, \,  a_1>a_2>a_3 \in \ZZ^3$$ Si $k_1\geq k_2 \geq k_3
\in \ZZ^3$, On note $\pi_{1,\infty}(k_1,k_2,k_3)$ la
représentation de $\U(3)(\RR)$ sur $V_{k_1,k_2,k_3}(\CC)$ (\ref{notations}), déduite de
l'inclusion $\U(3)(\RR) \subset \U(3)(\CC)\simeq \GL_3(\CC)$ 
donnée par $E \subset \CC$. Son
$L$-paramètre a pour changement de base à $\C^*$ le morphisme plus haut avec
$(a_1,a_2,a_3):=(k_1+1,k_2,k_2-1)$ (\cite{Rog1} \S 3.1).

\subsection{Existence de $\pi(\chi_0)$} \label{representationassocieeachi}

\begin{prop} \label{endo-1} Supposons que $k>1$, et que
$\varepsilon(\chi_0,1/2)=-1$,
alors il existe une unique représentation automorphe $\pi(\chi_0)$ de $\U(3)$ 
dont le changement
de base à $E$ a pour $L$-paramètre le morphisme $W_E \rightarrow \GL_3(\CC)$
défini par
	$$ \left( \begin{array}{ccc} \chi_0 |.|^{-1/2} &
& \\ & 1 & \\ & & \chi_0 |.|^{1/2} \end{array} \right ) $$
et telle que \begin{itemize}
\item Pour toute place $v, \pi(\chi_0)_v=\pi_v^n(\chi_0)$, 
\item $\pi(\chi_0)_{\infty}=\pi_{1,\infty}(\frac{k-1}{2},\frac{k-1}{2},1)$.
\end{itemize}
\end{prop}

{\it Preuve:} D'après \cite{Rog2}, page 397, le $A$-paquet 
$\Pi(\chi,1)$ existe pour $\U(3)$, car  $k > 1$.
D'après \cite{Rog2}, page 395 et 396, les $A$-paquets locaux 
correspondants sont les singletons $\{\pi_p^n(\chi_0)\}$ quand $p$ est 
décomposé (cf.~\ref{descpidec}) et des paires
$\{\pi_p^n(\chi_0),\pi_p^s(\chi_0)\}$ quand $p$ est inerte ou
ramifié (où $\pi_p^n$ a été définie en~\ref{descpiinerte}).
Enfin, d'après   \cite{Rog2} page 397, le $A$-paquet local à l'infini
 est un singleton $\{\pi(\chi_0)_{\infty}\}$ avec 
$\pi(\chi_0)_{\infty}:=\pi_{1,\infty}(\frac{k-1}{2},\frac{k-1}{2},1)$.

Considérons la représentation de $\U(3)(\A)$ :
$$\pi(\chi_0)=(\bigotimes_{p \text{ premier}} \pi_p^n(\chi_0) )\otimes 
\pi(\chi_0)_\infty.$$
D'après \cite{Rog2}, théorème 1.2, la multiplicité de $\pi$  dans le
spectre automorphe de $\U(3)$ est $(1+\varepsilon(\chi_0,1/2)(-1)^N)/2$, 
où $N$, défini
dans \cite{Roglivre} page 243, est le nombre de places à l'infini de
$\Q$ où $\U(3)$ est compact; on a donc $N=1$, et la multiplicité de 
$\pi(\chi_0)$ est donc 1.   
$\square$ \par \vspace{1 mm}

\subsection{Types}
  
\begin{prop}\label{types} Pour $p$ un nombre premier divisant $disc(\chi_0)$, il existe un groupe
$K_J=K_J(p)$ de $\U(3)(\Q_p)$, et une représentation irréductible
$J=J(p)$  de $K_J$ tels que
\begin{itemize}
\item $\Hom_{K_J}(J,\pi_p^n(\chi_0)\otimes(\chi_0^{-1}\circ\det)) \not = 0 $.
\item Pour toute représentation lisse irréductible $\pi$ de
$\U(3)(\Q_p)$ vérifiant $\Hom_{K_J}(J,\pi) \not = 0$, pour toute
place $v$ de $E$ au-dessus de $p$, notant $\pi_{E_v}$ le changement de
base de $\pi$ à $E_v$, il existe trois caractères lisses 
non ramifiés $\phi_1, \, \phi_2,\, \phi_3:
{E_v}^* \rightarrow \CC^*$, tels que
 $$L(\pi_{E,v})=\phi_1 \oplus \phi_2 \oplus  
\phi_3 \chi_0^{-1} .$$
En particulier, si $l \neq p$, $I':= \ker (\chi_l)_{|I_v}$, et si $\pi$ est la
composante en $p$ d'une représentation automorphe irréductible $\pi'$ de
$\U(3)$, alors $\rho_l(\pi')(I')=1$. 
\end{itemize}
\end{prop}
{\it Preuve: }
 Supposons d'abord que $p=v_1v_2$ est décomposé dans $E$. On note 
$G=\Gl_3(\Q_p)=\U(3)(\Q_p)$ (l'isomorphisme dépendant de la place $v_1$,
comme en~\ref{pdeccan}), $P$ le parabolique standard de type $(2,1)$, 
$M=\Gl_2(\Q_p)\times\Gl_1(\Q_p)$ son Levi (comme en~\ref{descpidec}), $N$ son 
radical unipotent, $\chi_0=\chi_{0,v_1}$ et $m$ la 
$p$-valuation du conducteur de $\chi_0$.

Notons $K_J$ le sous-groupe de $G$ des matrices dont la réduction 
modulo $p^m$ est de la forme
$$ \left(\begin{array}{ccc} * & * & * \\ * & * & * \\ 0 & 0 & y \end{array}\right), $$ et notons $J$ le caractère
complexe lisse de ce groupe qui à une matrice comme 
ci-dessus associe $\chi_{0}(y)^{-1}$. \par 

Par définition (cf. \ref{descpidec} pour la définition de $\lambda_p$), 
$\pi_p^n(\chi_0) \otimes \chi_0^{-1}\circ \det$ est la représentation
de $G$ sur l'espace  

$$V:=\{f :\ G\rightarrow \C,\ f\text{ lisse},\ \forall b \in P,\forall g \in G,
\, \, \ f(b g) = \lambda_p(b) 
\delta^{1/2}_P(b) \chi_0^{-1}(\det(b))
f(g)\}$$ donnée par $(g.f)(x)=f(x g)$. 
Définissons $f : G \rightarrow \CC$ par \begin{eqnarray*} 
f(bk) & = & \lambda_p(b) \delta^{1/2}_P(b) 
\chi_0^{-1}(\det(b)) \chi_0(k) \, \, \ \forall b \in P,\ k \in K_J,\\
f(g) &= & 0  \, \, \, \, \, \forall g \in G \backslash PK_J 
\end{eqnarray*}  
	On vérifie aisément que $f$ est bien définie, et qu'elle définit un
élément non nul de $\Hom_{K_J}(J,\pi_p^n(\chi_0)\otimes(\chi_0^{-1}\circ\det)).$
\par \vspace{1 mm}

Inversement, soit $\pi$ une représentation lisse irréductible de $G$, et supposons
que $$\Hom_{K_J}(J,\pi) \not=0.$$ 

Notons $B_J \subset K_J$ l'ensemble des matrices de $K_J$ qui sont 
triangulaires supérieures modulo $p$. On a alors $B_J \cap M = B \times \Gl_1(\Q_p)$ ou $B$ est le sous-groupe d'Iwahori standard (constitué des matrices triangulaires supérieures modulo $p$) de $\Gl_2(\Q_p)$. 
La restriction de $J$ à $B_J \cap M$ 
étant simplement le caractère $\chi_0^{-1}$ sur le second facteur 
$\Gl_1(\Q_p)$, le couple $(B_J \cap M,J_{|B_J \cap M})$ est un $\sss_M$-type
de $B_J$, où $\sss_M \in \BB(M)$ (le spectre de Bernstein de $M$, cf. 
\cite[page 772]{bus}) est la classe d'équivalence inertielle de $(T,\chi_0^{-1})$, $T$ désigne le tore maximal standard
de $\GL_3(\QQ_p)$ et $\chi_0^{-1}$ le caractère de $T$
envoyant $\diag(x,y,z)$ sur $\chi_0^{-1}(z)$. \par 

Il résulte immédiatement de la définition de {\it recouvrement (cover)} 
(\cite[8.1]{bk1}) que $(B_J,J_{|B_J})$ est un $G$-recouvrement de $(B_J \cap M, J_{|B_J \cap M})$. Le corollaire \cite[8.4]{bk1} (ou bien~\cite[page 55]{bk2})
assure alors que $(B_J,J_{|B_J})$ est un $\sss$-type pour $G$, où $\sss \in \BB(G)$ est la classe d'équivalence inertielle de $(T,\chi_0^{-1})$.
Comme $\pi$ contient $(K_J,J)$, elle contient aussi $(B_J,J)$, est son
 support cuspidale est donc de la forme $(\phi_1,\phi_2,\chi_0^{-1}\phi_3)$,
où les $\phi_i$ sont des caractères lisses non ramifiés de $\Q_p^\ast$.
Autrement dit, d'après les propriétés de la correspondance de Langlands 
locale, $L(\pi)$ a pour semi-simplification $\phi_1 \oplus \phi_2 
\oplus \chi_0^{-1}\phi_3$. Nous voulons maintenant montrer que $L(\pi)$ est 
semi-simple. Comme $\chi_0^{-1} \phi_3$ est ramifié, mais pas $\phi_1$ et $\phi_2$, 
on peut en tous cas écrire $L(\pi)=r \oplus \chi_0^{-1}\phi_3$, où $r^{ss}=\chi_1 \oplus \chi_2$, et nous voulons montrer que $r$ est semi-simple.

D'après ce qui précède et les propriétés de $L$, il existe une représentation
lisse irréductible $(\rho,W)$ de $\Gl_2(\Q_p)$ de support cuspidal $(\phi_1,\phi_2)$
telle que
  $$\pi\simeq 
\Ind_P(\rho \otimes (\phi_3 \chi_0^{-1})), \, \, \, L(\rho)=r$$
où  $\rho \otimes (\phi_3 \chi_0^{-1}$) désigne la représentation du 
parabolique standard sur $W$ qui à $(x,y)\in \Gl_2(\Q_p)\times\Gl_1(\Q_p)=M$ 
associe  $\phi_3(y) \chi_0(y)^{-1}\rho(x) \in\End(W)$. 
L'espace de cette induite est l'ensemble des fonctions $f:G \rightarrow W$, 
lisses, vérifiant $\forall m=(x,y) \in \Gl_2(\Q_p)\times\Gl_1(\Q_p)=M, \forall u \in N ,\forall g \in G,$
\begin{eqnarray} \label{equa1}
f(mug)=\delta^{1/2}_P(m) \chi_0^{-1}(y) \phi_3(y) \rho(x) f(g).
\end{eqnarray} 
Comme $\Hom_{K_J}(J,\pi) \not = 0$, il existe une fonction non nulle
$f:\ G \rightarrow W$ dans l'espace de 
$\Ind_P(\rho \otimes (\phi_3 \chi_0^{-1}))$ qui vérifie 
\begin{eqnarray}
\label{equa2}
(k.f)(g)=f(gk)=f(g)J_K(k)=f(g)\chi_0(k)^{-1}.
\end{eqnarray}
Choisissons $v=f(g) \in W$ non nul. 
 Combinant les équations~\ref{equa1} et \ref{equa2}, il vient 
$$\forall x \in \Gl_2(\Z_p), \rho(x) v = v.$$
La représentation $\rho$ est donc non ramifiée si bien que $r=L(\rho)$ est semi-simple, et finalement $L(\pi)=\phi_1 \oplus \phi_2 \oplus \phi_3 
\chi_0^{-1}.$ 

Le ``en particulier'' découle alors de la proposition~\ref{pastropram}

Le cas où $p$ est inerte ou ramifié dans $E$ est plus simple,
grâce aux travaux de L.Blasko (\cite{blas}) : $\pi_p^n(\chi_0)$ appartient à
la série principale, et son type défini en \cite[partie 7, page 181]{blas}
fait l'affaire. $\square$

\section{$I$-invariants et algèbre d'Atkin-Lehner}

\subsection{Notations} \label{notationsiwahori}
\subsubsection{} Dans toute cette section, $p$ est un nombre premier
fixé. $T$ désigne le tore diagonale de $G:=\GL_n(\QQ_p)$, $B$ son Borel supérieur,
$K:=\GL_n(\ZZ_p)$, $W \simeq \got{S}_n \subset K$, $I \subset K$ le sous-groupe d'Iwahori composé des éléments triangulaires supérieurs modulo $p$, $\Delta$ le
sous-groupe de $T$ des éléments à coefficients dans $p^{\ZZ}$, $\Delta^+$ le
sous-monoïde de $\Delta$ des éléments de la forme
$\diag(p^{a_1},\cdots,p^{a_n})$ avec $a_1\leq \cdots \leq a_n \in \ZZ$. 
Si $X$ est un groupe topologique localement compact, $\delta_X$ désigne le caractère module de $X$. 

\subsubsection{} Si $U$ est un sous-groupe compact ouvert de $G$, l'algèbre
de Hecke de $G$ relativement à $U$, $\mathcal{C}^{\infty}(U\backslash G/U)$, est l'algèbre de convolution des
fonctions complexes lisses à support compact sur $G$ invariantes à droite et à gauche
par $U$. Si $g \in G$, on note $[UgU] \in \mathcal{C}^{\infty}(U\backslash
G/U)$ la fonction caractéristique de de $UgU \subset G$, on prend la
convention que $[U]^2=[U]$. $\mathcal{C}^{\infty}(I\backslash G/I)$
est l'algèbre de Hecke-Iwahori, on note $\mathcal{A}(p)$ son sous-anneau de
engendré $\ZZ[1/p]$ et les fonctions caractéristiques $[IuI]$, $u \in \Delta^{+}$. Il est connu
(par exemple \cite{Rog3} \S 1) que $\mathcal{A}(p)$ est commutative, et que pour chaque $u,u' \in \Delta^{+}$,
$[IuI]$ est inversible dans $\mathcal{A}(p)$, $[IuI][Iu'I]=[Iuu'I]$ et
$IuIu'I=Iuu'I$. On pose

$$ u_i:=\diag(1,\dots,1,p,\dots,p) \in \Delta^{+}, \textrm{ où $p$ apparaît $i$ fois}, \, \,
0\leq i \leq n $$

\subsection{$I$-invariants des représentations non ramifiées} 

\subsubsection{} \label{iwa1} Soit
$\psi=(\psi_1,...,\psi_n): (\QQ_p^*)^n \rightarrow \CC^*$ un caractère lisse
et non ramifié,
$X(\psi)$ la représentation complexe lisse irréductible
non ramifiée de $G$ de $L$-paramètre valant sur le Frobénius géométrique:

\begin{eqnarray*} \left( \begin{array}{ccc} \psi_1(p)
& &  \\ & \cdots &  \\ & & \psi_n(p) \end{array} \right)
\end{eqnarray*}

Voyant $\psi$ comme un caractère complexe lisse de $B$ trivial sur les
unipotents supérieurs, on pose
$$\Ind_B(\psi):=\{f: G \rightarrow \CC, \textrm{ lisses},
f(bg)=\delta_B^{1/2}(b)\psi(b)f(g)\, \,  \forall b \in B\}$$ vue comme
représentation lisse de $G$ par translation à droite les fonctions. On sait
alors
 que $X(\psi)$ est l'unique sous-quotient irréductible non ramifié de $\Ind_B(\psi)$. Le
caractère de l'algèbre de Hecke non-ramifiée de
$G$ sur la droite des $K$-invariants de $X(\psi)$ (ou encore de
$\Ind_B(\psi)$, ce qui est la même chose) est bien connu: si $t_{p,i}$ est
la valeur propre de l'opérateur de Hecke $T_{p,i}:=[Ku_iK]|\det(u_i)|$ comme plus haut, on a

        $$ \prod_{i=0}^n (1-\psi_i(p)T)=\sum_{i=0}^n (-1)^i p^{i(i-1)/2} t_{p,i}
T^i$$

\subsubsection{} \label{iwa2} Nous allons commencer par décrire la représentation de
$\mathcal{A}(p)$ sur les $I$-invariants de $\Ind_B(\psi)$. On
verra les caractères de $T$ par restriction comme des caractères de
$\Delta^{+}$, puis de $\mathcal{A}(p)$. Pour tout caractère $\theta: T
\rightarrow \CC^*$, $\sigma \in W$, on dispose d'un caractère
$\theta^\sigma$ défini par $\theta^{\sigma}(t)=\theta(\sigma^{-1}t\sigma)$; ainsi $(\psi_1,...,\psi_n)^\sigma=(\psi_{\sigma(1)},...,\psi_{\sigma(n)})$. \par \vspace{2 mm}

\begin{lemme} \label{U} (\cite{Che} 4.8.4) La semi-simplification de
$\Ind_B(\psi)^I$
comme $\mathcal{A}(p)$-module est $$\oplus_{\sigma \in W}
\delta_B^{1/2}\psi^{\sigma}$$ \end{lemme}

{\it Remarques:} Le calcul ci-dessus et la formule pour le polynôme de
Hecke de $X$ montrent que si $\sigma \in W$, $\psi^{\sigma}(U_i)$ (noter la
disparition du $\delta_B^{1/2}$) est un produit de $i$ valeurs propres
``distinctes" du Frobénius géométrique dans la représentation de $W_{\QQ_p}$
attachée à $X$ par la correspondance de Langlands locale non ramifiée. \par \vspace{3 mm}

        En particulier, supposons $\Ind_B(\psi)$ irréductible, alors
 $\Ind_B(\psi)=X(\psi)$, $\dim_{\CC}(X(\psi)^I)=n!$ et l'action de
$\mathcal{A}(p)$ sur les $I$-invariants de $X$ est calculée par le
lemme. Ceci se produit en particulier quand $\psi$ est
essentiellement tempérée (i.e $|\psi_i(p)|$ indépendant de $i$, \cite{Z} 4.2) 
mais pas pour la représentation $\pi_p^n(\chi_0)$ introduite en \S
\ref{descpidec} ($p$ décomposé dans $E$ et ne
divisant pas $\cond(\chi_0)$). Pour traiter ce cas là, nous aurons
besoin du résultat général suivant, impliquant par ailleurs le
lemme précédent. 
\subsubsection{} Soit $P=MN \subset B$ un parabolique de Levi $M$,
$\psi: M \rightarrow \CC^*$ un caractère non ramifié,
$\Ind_P(\psi)$ l'induite parabolique lisse normalisée:

$$\Ind_P(\psi):=\{f: G \rightarrow \CC, \textrm{lisses},
f(pg)=\delta_P^{1/2}(p)\psi(p)f(g)\, \,  \forall p \in B\}$$

\par \vspace{1 mm}

 Soit $W_P \subset W$ le sous-groupe de Coxeter correspondant à $P \subset
G$, pour chaque $\sigma \in W$ on choisit l'unique élément dans $W_P.\sigma$ de
longueur minimale (\cite{Hum} prop. 1.10 c)), et on note $W(P)$ l'ensemble des représentants
de $W_P\backslash W$ obtenu. On pose $z:=\delta_P/\delta_B=\delta_{B\cap
M}$.

\begin{lemme} \label{Ubis} La semi-simplification de $\Ind_P(\psi)^I$
comme $\mathcal{A}(p)$-module est $$\oplus_{\sigma \in W(P)}
\delta_B^{1/2}(z^{1/2}\psi)^\sigma$$
\end{lemme}
{\it Remarques:} Le lemme montre que l'unique sous-quotient non ramifié de $\Ind_P(\psi)$ est isomorphe
à $X(\psi . z^{1/2})$. Si $\Ind_P(\psi)$ est irréductible (voir \cite{Z} 3.2,4.2), 
il coïncide donc avec $X(\psi.z^{1/2})$ dans ce cas. \par
\vspace{2 mm}

{\it Preuve:} On pose $X'=\Ind_P(\psi)^I$, c'est un module sous
l'algèbre de Hecke-Iwahori de $G$, en particulier sous
$\mathcal{A}(p)$. La décomposition de Bruhat-Iwahori $G=\coprod_{\sigma
\in W(P)} P\sigma I$ montre que $\dim_{\CC}(X')=|W|/|W_P|$. On considère
la $\CC$-base de $X'$ suivante: si $\sigma \in W(P)$, $e_\sigma$ est
l'élément de $X'$ nul hors de $P\, \sigma \,I$ et tel que $e_\sigma(\sigma)=1$.
$e_\sigma(\sigma')=0$ ou $1$ selon que $\sigma' \in W_P.\sigma$ ou non . Si $\sigma' \in
W$, on commence par calculer $[I\sigma'I](e_1)$. La décomposition de
Bruhat-Iwahori, ainsi que la multiplication des cellules, montre
que $(I\sigma'I) \cap (\sigma^{-1}PI) = \emptyset $ à moins que $\sigma \in
W_P{\sigma'}^{-1}$. En particulier, $[I\sigma'I](e_1)=a_{\sigma'}e_{{\sigma'}^{-1}}$
où $a_{\sigma'} \in \C^*$ (car $[I\sigma'I]$ est inversible dans l'algèbre de Hecke-Iwahori).
\par

La décomposition de Bruhat-Iwahori montre que si $u \in \Delta^{+}$,
$(IuI)\cap \sigma^{-1}PI$ est vide à moins que $\sigma \in W_P$, ce qui
implique que $e_1$ est propre sous l'action de $\mathcal{A}(p)$,
de caractère $[IuI] \mapsto (\delta^{1/2}\psi)(u)|(IuI\cap
PI)/I|$. Quand $u_i=u$, on peut calculer $|(IuI\cap
PI)/I|=|(I^u\cap PI)/(I^u\cap I)|$, $I^u:=u^{-1}Iu$, on vérifie
qu'il vaut $z^{1/2}(u)$.
\par Par les relations de Bernstein (\cite{Rog3} \S 1, \S 5), on en
déduit que si l'on ordonne la base des $e_\sigma$, $\sigma \in W(P)$, par
ordre croissant avec la longueur de $\sigma$, l'action de
$\mathcal{A}(p)$ est triangulaire supérieure. Toujours par les
relations de Bernstein, on trouve alors les caractères de
l'énoncé (voir par exemple \cite{Che} \S 4.8.4 pour des détails
supplémentaires). $\square$

\par \vspace{2 mm}

\subsubsection{} \label{exemple}  L'exemple qui nous intéresse est le cas de la
représentation non ramifiée $X(\psi):=\pi_p^n(\chi_0)$ définie en \ref{descpidec}, 
$\pi_p^n(\chi_0)$ est l'induite du parabolique standard $P$ de type $(2,1)$ du
caractère $\chi_0(\det(.)) \times 1$. On trouve $W(P)=\{1,(3, 2), (3, 2,
1)\}$. Si $\sigma \in W(P)$, le triplet
$(\psi^{\sigma}(u_1),\psi^{\sigma}(u_2/u_1),\psi^{\sigma}(u_3/u_2))$ associé
à $\sigma$ est alors explicitement donné par:
\begin{center} $\sigma=1$, $(1,\chi_{v_1}^{\bot}(p),\chi_{v_1}(p))$ \end{center}
\begin{center} $\sigma=(3, 2)$, $(\chi_{v_1}^{\bot}(p),1,\chi_{v_1}(p))$ \end{center}
\begin{center} $\sigma=(3, 2, 1)$, $(\chi_{v_1}^{\bot}(p),\chi_{v_1}(p),1)$ \end{center}

\section{Déformations des représentations cristallines raffinées}

\subsection{Raffinement d'une représentation cristalline} 
\subsubsection{} Soit $F$ un corps local, $V$ un $F$-espace vectoriel de
dimension finie muni d'une représentation continue de $\Dp$. Nous supposerons que $V$
est cristalline, que ses poids de Hodge-Tate $k_1 < ... < k_n$ sont tous distincts
et que les valeurs propres du Frobénius de $D_{cris}(V)$ sont dans $F$ (cf.
\ref{notationshodge} pour les conventions). Imitant Mazur (\cite{Ma}), on appellera {\it raffinement} de $V$ la donnée d'un
ordre $\mathcal{R}:=(\varphi_1,...,\varphi_n) \in F^n$ sur les valeurs propres du
Frobénius de $D_{cris}(V)$. On notera $(V,\mathcal{R})$ la représentation
$V$ munie de son raffinement $\mathcal{R}$. \par \vspace{2 mm} 

\subsubsection{} \label{raff2}La donnée d'un raffinement $\mathcal{R}$ de $V$ nous permet de définir des
$F_i(\mathcal{R}):=\varphi_i/p^{k_i} \in F^*$ et des $U_i(\mathcal{R}):=\prod_{j=1}^iF_j$, $1 \leq i \leq n$.  La donnée de
tous les $F_i(\mathcal{R})$,
ou encore celle des $U_i(\mathcal{R})$, est bien sur équivalente à celle de $\mathcal{R}$.
Notons que $F_i(\mathcal{R})$ est encore une valeur propre du Frobénius de $D_{cris}(V(k_i))$, et
$U_i(\mathcal{R})$ en est une de celui de $\Lambda^i(D_{cris}(V))
\simeq D_{cris}(\Lambda^i(V))$. On rappelle que la formation de $D_{cris}$ commute aux opérations
tensorielles sur les représentations cristallines (\cite{Fo2} 1.5.2, 5.1.2).  \par
\vspace{2 mm}

\subsubsection{} {\it Remarques:} i) Si le polynôme caractéristique du Frobénius 
de $D_{cris}(V)$ a $n=\dim_F(V)$ racines distinctes, alors $V$ admet exactement $n!$
raffinements. C'est la cas par exemple si $V$ est ordinaire (cf. \cite{PR}). Dans ce cas, on
dispose de plus d'un raffinement canonique donné par $|\varphi_i|=p^{-k_i}$,
appelé raffinement "ordinaire" (il ne nous sera pas utile dans la suite).  \par \vspace{1 mm}

ii) Un des intérêts essentiels de la notion de raffinement dans ce texte 
vient de ce que la théorie des familles $p$-adiques de formes
automorphes produit des déformations de représentations cristallines
raffinées de $\Dp$. Il faut bien noter que de telles déformations d'une même représentation $V$
mais partant de raffinements distincts sont en général très différentes (par
exemple l'une peut être génériquement irréductible, l'autre non). 

\subsection{Raffinements et algèbre d'Atkin-Lehner} 
\subsubsection{} \label{raff-notations} Par commodité
d'exposition, nous nous restreignons à $\U(3)$ plutôt qu'à un groupe unitaire
quelconque, c'est de toutes façons suffisant pour les objectifs de ce texte. Soit $p$ un nombre premier décomposé dans $E$, $\iota: \CC \rightarrow
\overline{\QQ}_p$ un isomorphisme
de corps, $v_1$ la place de $E$ au dessus de $p$ donnée par $E \subset \CC
\overset{\iota}{\rightarrow}\overline{\QQ}_p$, $v_2$ l'autre place, $p=v_1v_2$. La donnée de $v_1$
nous permet de plus d'identifier $\U(3)(\QQ_p)$ à $\GL_3(\QQ_p)$ comme en
\ref{descpidec}. Fixons $\Pi$ une
représentation automorphe irréductible de $\U(3)$ telle que $\Pi_p$
est non ramifiée et $\Pi_{\infty}=\pi(k_1\geq k_2 \geq k_3)$. Comme dans
\ref{iwa1}, $\Pi_p=X(\psi)$ pour un certain caractère non ramifié $\psi:
T \rightarrow \CC^*$. \par \vspace{2 mm}
\subsubsection{} \label{raff-galoisienne} On a vu dans 
\ref{propgal} que la donnée de $\iota$ permet d'associer à $\Pi$ une
représentation semi-simple continue $\rho_p: \Gal \rightarrow
\GL_3(F)=\GL_F(V)$, $V$ étant un espace
vectoriel de dimension $3$ sur un certain corps local $F$. On rappelle que 
$D_{v_i}$ est un groupe de
décomposition dans $\Gal$ associé à la place $v_i$ de $E$, et on note
$V_i$ la représentation continue de $\Dp$ sur $V$ obtenue par restriction de
$V$ à $D_{v_i}$. Puisque $\rho_p^{\bot} \simeq \rho_p$, on sait que $V_2 \simeq V_1^{\bot}$. De plus, par
la proposition \ref{nonramcris} \S \ref{classification}, $V_1$ est
cristalline de poids de Hodge-Tate $-k_1-1 < -k_2 < -k_3 +1$, et son
Frobénius cristallin a même polynôme caractéristique (modulo $\iota$) que
l'image du Frobénius géométrique de $W_{\QQ_p}^{nr}$ dans le $L$-paramètre
de $X(\psi)$, i.e $\prod_{i=1}^3(X-\iota(\psi_i(p)))$. Nous allons
donner une interprétation automorphe de certains raffinements de $V_1$ en
terme de $\Pi_p$. \par \vspace{2 mm} 
\subsubsection{} \label{accessible} La semi-simplification de $\mathcal{A}(p)$ agissant sur $\Pi_p^I$ a
été calculée en \ref{iwa2}, c'est une somme de caractères de
$\mathcal{A}(p)$ de la forme $\delta_B^{1/2}\psi^\sigma$, pour certains $\sigma \in
W$. On dira que $\sigma$ est {\it accessible} pour $\Pi$ si $\delta_B^{1/2}\psi^\sigma$ apparaît, cela ne dépend
que de $\Pi_p$, et il est équivalent de demander qu'il existe 
un vecteur $v \in \Pi^I$ sur lequel $\mathcal{A}(p)$ agisse par $\delta_B^{1/2} \psi^\sigma$. 
Pour tout $\sigma \in W$, le lemme \ref{U} \S \ref{iwa2} montre que l'on construit un raffinement
$\mathcal{R}(\sigma)$ de $V_1$ en posant 
$$\mathcal{R}(\sigma):=(\psi_{\sigma(3)}(p),\psi_{\sigma(2)}(p),\psi_{\sigma(1)}(p))$$

Un raffinement de $V_1$ sera dit {\it accessible} s'il est de la forme 
$\mathcal{R}(\sigma)$ avec $\sigma$ accessible pour $\Pi$. Si $\Pi_p$ est irréductible, on a vu
en \ref{U} que tous les raffinements de $V_1$ sont
alors accessibles. Cela se produit en particulier quand $\Pi_p$ est
tempérée. Par contre, si $\Pi=\pi(\chi_0)$, les raffinements accessibles sont
ceux donnés dans \ref{exemple}. Dans ce cas précis, on remarque par exemple
que $V_1$ est ordinaire mais que le raffinement ordinaire n'est pas accessible. \par \vspace{2 mm}

\subsubsection{} \label{versionU} Soit $s=(k_1\geq k_2 \geq k_3) \in \ZZ^3$, notons $\nu_s$ le
caractère $\Delta \rightarrow p^{\ZZ}$ sur le vecteur de plus haut poids de
$V_{s}^*(\QQ)=V_{-s}(\QQ)$. D'après \ref{raff-galoisienne} et
\ref{plushautpoids}, si $1 \leq i \leq 3$, 
$(\delta_B^{-1/2}\nu_s)(u_i/u_{i-1})$ est une puissance de $p$ d'exposant le $i^{ieme}$ poids de Hodge-Tate
(rangés par ordre croissant) de $V_1$. On pose 

$$U_i^s:=\frac{[Iu_iI]}{\nu_s(u_i)} \in \mathcal{A}(p)$$

Avec ces notations, $\mathcal{R}(\sigma)$ est encore le raffinement de $V_1$ défini
par ses $U_i(\mathcal{R}(\sigma))$ avec la formule $U_i(\mathcal{R}(\sigma)):=(\delta_B^{1/2}\psi^\sigma)(U_i^s)$. 

\subsection{Variation de représentations cristallines raffinées, d'après M.Kisin} 

\subsubsection{} \label{notationkisin} Soit $F \subset \CC_p$ un corps
local, $X$ un $F$-affinoide réduit d'anneau $A(X)$, la norme réduite de $A(X)$ en fait une
algèbre de Banach. Soit $M:=A(X)^n$, $G$ un groupe topologique, $\rho: G \longrightarrow \GL(M)$ une représentation continue. 
Si $F \subset E \subset \CC_p$ est un sous-corps complet, $x \in X(E)$, on note 
$M_x:=M \otimes_{A(X)} E$, $A(X) \rightarrow E$ étant l'évaluation en $x$, et $\rho_x: G \rightarrow
\GL_n(E)=\GL_E(M_x)$. \par \vspace{1 mm}

\subsubsection{} \label{generalitesfamilleslocales} 
Soit $X$ un $F$-affinoide réduit d'anneau $A(X)$, on se donne \par \vspace{ 2
mm}\begin{itemize} 
\item  a) $\kappa=(\kappa_1,...,\kappa_n): X \longrightarrow \A^n$, 
un $F$-morphisme analytique, 
\par \vspace{1 mm}

\item b) $Z \subset X(F)$ un sous-ensemble Zariski-dense tel que $\kappa(Z) \subset
\ZZ^{3,--}$. \par \vspace{1 mm}

\item c) $M:=A(X)^n$, et  $\rho: \Dp \longrightarrow
\GL_{A(X)}(M)$ une représentation continue, \par \vspace{1 mm}

\item d) $F_1,...,F_n \in A(X)$, \par \vspace{1 mm}
 
\end{itemize}

On suppose de plus que: \par \vspace{2 mm}

i) Pour tout réel $C >0$, $\{ z \in Z, \,\, -\delta(\kappa(z))> C \}$ est
Zariski-dense dans $X(F)$, \par
\vspace{1 mm}
ii) Pour tout $x \in X(F)$, $\kappa(x) \in F^n$ est l'ensemble des poids de
Hodge-Tate-Sen de $M_z$, rangés par ordre croissant si $x \in Z$, \par \vspace{1 mm}
iii) Pour tout $z \in Z$, $M_z$ est cristalline, \par \vspace{1 mm}
iv) Si $z \in Z$, $(p^{\kappa_1(z)}F_1(z),...,p^{\kappa_n(z)}F_n(z))$ est 
un raffinement de $M_z$, \par \vspace{1 mm}
v) $|F_i|$ est constant sur $X$, \par \vspace{1 mm}

vi) Deux éléments de $\kappa(Z)$ diffèrent d'un élément de $(p-1)\ZZ^n$.
\vspace {3 mm}

Faisons quelques remarques sur ces hypothèses. L'existence de
$\kappa$ satisfaisant $ii)$ n'est pas directement conséquence des travaux de
Sen (\cite{Sen1}, \cite{Sen2}), à cause de l'hypothèse "rangés par
ordre croissant sur $Z$". Nous ne savons pas dans quelle mesure l'existence de $F_i$ satisfaisant
iv) est automatique (condition sur $\kappa$ ?), cf. \cite{Ki} à ce sujet. 
Nous ne savons pas non plus si vi) est automatique. Dans la pratique (cf.
\S \ref{defo}), ces hypothèses seront satisfaites et $\kappa$ sera un morphisme fini sur son image, dominant restreint à chaque composante
irréductible de $X$. \par \vspace{1 mm}

\subsubsection{} On se place sous les hypothèses a), b),
c), ii), iii) et iv) de \ref{generalitesfamilleslocales}. On notera $A(X)(\kappa_i)$ le $A(X)$-module $A(X)$ 
sur lequel on fait agir continument $\Dp$ par le composé du  caractère cyclotomique par
caractère suivant de $\ZZ_p^*$: $$x \rightarrow (x/\tau(x))^{\kappa_i} \tau(x)^{n_i} \in
A(X)^*$$ où  $n_i=\kappa_i(z)$ est un entier bien défini modulo $(p-1)\ZZ$ indépendamment
de $z \in Z$ par l'hypothèse $vi)$, $\tau: \ZZ_p^* \rightarrow
\mu_{p-1}(\QQ_p)$ la réduction modulo $p$ composée par le caractère de
Teichmüller. Si $z\in Z$, $A(X)(\kappa_i)_z$ est un caractère cristallin
de $\Dp$. \par \vspace{1 mm}

Si $N$ est un $A(X)$-module de libre type fini muni d'une représentation
continue de $\Dp$, on notera $N(\kappa_i):=N \otimes_{A(X)} A(X)(\kappa_i)$ vue comme
représentation de $\Dp$. On a donc défini les $M(\kappa_i)$, leurs
évaluations en $z \in Z$ sont cristallines de poids
$(\kappa_1(z)-\kappa_i(z),...,\kappa_n(z)-\kappa_i(z)) \in \ZZ^{n,--}$. \par \vspace{2 mm}

\subsubsection{}On se place dans les hypothèses de
\ref{generalitesfamilleslocales}:

\begin{prop} \label{kisin} Soit $x \in X(F)$, $\kappa(x)=(k_1,...,k_n)$ alors
$$D_{cris}((\Lambda^i V_x)(k_1+...+k_i))^{\varphi=F_1(x)\cdots F_i(x)} \,
\, \, \, \textrm{ est non
nul.} $$ 
Si de plus les $k_i$ sont entiers, alors
$$D_{cris}(\Lambda^i(V_x))^{\varphi=\prod_{j=1}^ip^{k_j}F_j(x)} \, \, \, \,
\textrm{ est non         
nul} $$

\end{prop}

{\it Preuve:} La seconde assertion découle de la première.
Soit $i \in \{1,...,n\}$ fixé, on pose 
$N:=\Lambda^i(M)(\kappa_1)\cdots(\kappa_i)$. Rappelons que, comme dit plus
haut (\ref{raff2}), la formation de $D_{cris}(V)$ commute aux opérations tensorielles sur
la catégorie des représentations cristallines. En
particulier, si $x \in Z$, on en déduit que $N_x$ est cristalline, 
de plus petit poids de Hodge-Tate égal à $0$, et que $\prod_{j=1}F_j(x)$ est
une racine de son Frobénius cristallin. Rempla\c cant $\kappa$ par
$(\kappa_{J_1},...,\kappa_{J_r})$, les $J_j$ parcourant l'ensemble ordonné
lexicographiquement des $r={{n}\choose{i}}$ parties à $i$ éléments de
$\{1,...,n\}$, $\kappa_{(m_1 \leq ... \leq m_i)}:=\sum_{j=1}^i
\kappa_{m_j}$, les conditions $i)$, $ii)$, $iii)$, 
$iv)$ et $v)$ sont satisfaites pour $N$ à la place de $M$.

En fait, seule la première valeur propre va servir pour la suite. On
est ramené à prouver la proposition pour $M$ et $i=1$, et on peut supposer que le premier poids de Hodge-Tate
des points dans $Z$ est nul. \par 
        A partir de là, on suit Kisin \cite{Ki} 6.3. Par l'hypothèse $iv)$
et 5.14 il suffit donc de vérifier les hypothèses de sa proposition 5.13. Dans ses
notations $M:=M$, $I:=Z$, et si $z \in X(F)$ est d'idéal maximal $m_z$,
$R_z:=A/m_z=F$. Pour chaque entier positif $k$, on définit $I_k$ comme étant 
l'ensemble des points dans $Z$ ayant leur second poids de Hodge-Tate plus grand que $k$, et tels que
$v(F_1(x))<k$. $I_k$ est Zariski-dense dans $X(\CC_p)$ par $i)$, il
satisfait donc les conditions 2 et 3 de {\it loc. cit}. 5.13. Soit $x \in I_k$, par hypothèse
$M_x$ a pour plus petit poids de Hodge-Tate $0$ et $D_{cris}(M_x)^{\varphi=F_1(x)}$ est non nul,
il existe donc une application $\Dp$-équivariante non nulle
$M_x^* \rightarrow (\Fil^0(B_{cris})\otimes_{\QQ_p}F)^{\varphi=F_1(x)}$. 
Cette application ne se factorise pas par $Fil^k(B_{cris}) \otimes_{\QQ_p} F$ car
par faible-admissibilité de $D_{cris}(M_x)$, $Fil^k(D_{cris}(M_x))^{\varphi=F_1(x)}$ est encore faiblement admissible, 
de plus petit poids de Hodge-Tate $k$, mais $v(F_1(x))<k$ par hypothèse. On
conclut mot à mot comme dans les $3$ dernières lignes de la preuve de
{\it loc. cit}. 6.4. $\square$ \par \vspace{2 mm}

{\it Remarques: } i) En général, il est bien sur faux que sous les hypothèses
de la proposition \ref{kisin}, $D_{cris}(V(k_i))^{\varphi=F_i(x)}$ est non nul si 
$i>1$. La considération d'une famille $p$-adique de formes modulaires
passant par une forme modulaire parabolique propre de niveau $\Gamma_0(p)$ qui est
$p$-nouvelle en donne un contre-exemple. \par \vspace{1 mm}
ii) Ce que l'on a fait pour les $\Lambda^i$ vaut de même aussi
naturellement pour n'importe quel foncteur de Schur de $\GL(n)$, c'est en
fait une conséquence de la proposition ci-dessus. \par \vspace{1 mm}
iii) Si $n=3$, l'isomorphisme canonique $V^* \otimes \det(V) \simeq
\Lambda^2(V)$ permet de reformuler la proposition pour $i=2$ en terme de $V^*$.
Elle montre alors que $D_{cris}(V(k_1))^{\varphi=F_1(x)}$ et
$D_{cris}(V^*(-k_3))^{\varphi=F_3(x)^{-1}}$ sont non nuls.

\section{Extensions et pseudo-représentations}

\subsection{Existences de réseaux stables} \label{bidule}
\subsubsection{} \label{semisimples} Soit $K$ un corps de caractéristique
$0$, $V=K^r$, $G$ un groupe, $\rho: G  \rightarrow \GL(V)$ une représentation
semi-simple. Soit $(\rho=\oplus_{i=1}^n \rho_i,V=\oplus_{i=1}^n W_i)$ la
décomposition isotypique de $(\rho,V)$ : pour chaque $i$, $(W_i,\rho_i)$
est somme directe de $n_i$ copies d'une représentation irréductible 
$(V_i,\rho'_i)$, et $\rho'_i \not \simeq \rho'_j$ si $i \neq j$.

On dira que $V$ satisfait la condition (ABS) si les $\rho'_i$ sont {\it 
absolument} irréductibles.  \par \vspace{2
mm}
	On suppose dès maintenant et dans tout \S \ref{bidule} que $V$ satisfait (ABS).
On fixe de plus $A$ un sous-anneau intègre noethérien de $K$ tel que $K=\textrm{Frac}(A)$ et
$\tr(\rho(G)) \subset A$.

\subsubsection{} \label{reseaux} Si $B$ est un sous-anneau de $K$, on appelle 
$B$-réseau de $V$ un sous-$B$-module libre $\Lambda$ tel que $K.\Lambda =
V$. Si $V$ est une représentation de $G$, on dit 
qu'un $B$-réseau est stable s'il est stable sous l'action de $G$.

\begin{lemme} \label{stable} 
i) L'image de $A[G]$ dans $End_K(V)$ est de type fini sur $A$. \par \vspace{2 mm} 
ii) Si $A$ est normal, $\tr(\rho'_i(G)) \subset A$. \par \vspace{2 mm}
iii) Si $A$ est principal, $V_i$ admet un $A$-réseau stable. \par \vspace{2 mm} 
iv) Plus généralement, si $P \in \Spec(A)$ tel que $A_P$ est de valuation
discrète, il existe $g \in A\backslash P$ tel que $V_i$ admet un $A_g$-réseau stable. \par \vspace{2 mm}
v) Supposons de plus que $A$ est soit une algèbre affinoide, soit local
complet, $G$ un groupe topologique,  si $T: G \rightarrow A$ est continue, alors les $g \rightarrow \tr(\rho'_i(g))$ sont continus. De plus,
tout $A$-réseau stable est une représentation continue de $G$.
\end{lemme}

{\it Preuve:} Par la théorie des modules semi-simples, l'image de $K[G]$ dans $\End(V)$ est somme directe de ses images dans les
$\End(W_i)$. De plus, l'image de $K[G]$ dans $\End(W_i)$ est l'action
diagonale de $\End(V_i)$ dans $\End(W_i)$. En particulier, on dispose d'un
$e_i \in K[G]$ tel que $\rho(e_i)$ est le projecteur $G$-équivariant sur
$W_i$. Aussi, $$\exists f \in K^*, \, \, \forall g \in G, \,\,
\tr(\rho_i(g))=\tr(\rho(e_ig)) \in f A $$
On déduit la même assertion pour $\rho'_i$, en remplaçant $f$ par $f/n_i$.
\par \vspace{1 mm}
Soit $i$ fixé, $d=\dim(V_i)$, on fixe une base de $V_i$ nous permettant de
l'identifier à $K^d$. La représentation $\rho_i'$ étant absolument 
irréductible, un
théorème de Wedderburn assure l'existence de $g_1,...,g_{d^2} \in G$, tels
que les $\rho_i'(g_k)$ engendrent $\End(V_i)$ comme $K$-espace vectoriel.
Ainsi, $$M:=((\tr(\rho_i'(g_kg_l)))_{1\leq k,l \leq d^2}) \in \GL_{d^2}(K)$$
Soit $f' \in K^*$ tel que $M^{-1}$ et les $\rho_i'(g_k)$ soient à
coefficients dans $f'A$, alors pour tout $g \in G$ on a $\rho_i'(g) \in ff'/n_i M_d(A)$. 
Autrement dit, 
\begin{equation} \label{ligne} A^d \subset \rho_i'(G).A^d \subset ff'/n_i A^d 
\end{equation}
 En particulier, $A[\rho_i'(G)]$ s'injecte dans
$\textrm{Hom}_A(A^d,(ff'/n_i)A^d)$. Comme $A$ est noethérien, cela prouve $i)$.
De plus, si $A$ est principal, \eqref{ligne} montre que $A[\rho_i'(G)].A^d$ est un réseau
stable, cela montre $iii)$. On en déduit $ii)$ car un anneau normal est intersection
d'anneaux de valuation discrète. $iv)$ est une conséquence de $ii)$ et $i)$.
\par \vspace{1 mm}
La première assertion de $v)$ découle de $n_i\tr(\rho'_i(g))=\tr(\rho(e_ig))$,
prouvons la seconde. Si $\Lambda$ est un réseau stable de $V$, il suffit de vérifier que $\psi: G \rightarrow
\End_A(\Lambda)$ est continue, car alors $g \mapsto \psi(g)^{-1}=\psi(g^{-1})$ le sera aussi.
On rappelle que si $A$ est affinoide (\cite{BGR} 3.7.3) ou local
complet, tout $A$-module de type fini a une topologie canonique et 
que toute application $A$-linéaire entre deux tels modules est continue et fermée. 
Par $i)$, on peut trouver $g_1,...,g_s \in G$ engendrant $M:=A[\rho(G)]$, on
munit $M$ de la topologie discutée ci-dessus. Il suffit de montrer que 
$\psi^*: G
\rightarrow M$ est continue. Par semi-simplicité
de $V$ comme $G$-représentation, l'application $\psi^{**}:G \rightarrow A^s$, $g \rightarrow (\tr(g_ig))$, 
induit une injection $A$-linéaire $M \rightarrow A^s$, nécessairement
continue et fermée. La continuité de $\psi^{**}$ conclut. $\square$ \par
\vspace{3 mm}

\subsection{Représentations attachées aux pseudo-caractères}

\subsubsection{} \label{generalitespseudo} Soit $G$ un groupe, $A$ un
anneau commutatif, on rappelle qu'une fonction $T: G \longrightarrow A$ est un
pseudo-caractère sur $G$, de dimension $n \in \NN$, à coefficients dans $A$ si 

$$\forall g,h \in G, \, \, \, T(gh)=T(hg)$$
$$\forall g=(g_1,...,g_{n+1}) \in G^{n+1}, \, \, \sum_{\sigma=c_1...c_r \in \got{S}_{n+1}} \epsilon(\sigma) \prod_{i=1}^{r}
f(c_i(g))=0$$
$$ \textrm{ $n$ est le plus petit entier ayant la propriété ci-dessus} $$
Ici $\sigma=c_1...c_r$ est la décomposition en cycles de $\sigma$, et si
$c=(j_1,...,j_s)$ est un cycle $c(g)=\prod_{i=1}^s g_{j_i}$, cf. \cite{Tay}
\S 1, \cite{Rou} \S 2. \par \vspace{1 mm}

	La trace d'une représentation $\rho: G \rightarrow \GL_n(A)$ est en
particulier un pseudo-caractère sur $G$, à coefficients
dans $A$ (cf. {\it loc. cit.}); il est de dimension $n$ si $A$ est intègre de
caractéristique $0$ (cf. \cite{Rou} 2.4). On discute de l'assertion
réciproque dans les paragraphes qui suivent.
	
\subsubsection{} \label{repattacheesauxpseudorep} Si $A=F$ est un corps, il est connu que quitte à
faire une extension séparable finie de $F$, $T$ est la trace d'une 
représentation semi-simple $G \rightarrow \GL_n(F)$, unique à isomorphisme
près, satisfaisant la propriété
(ABS) de~\ref{semisimples} (voir \cite{Tay} \S 1
théorème 1, \cite{Rou} 4.2 dans cette généralité). On dira que $T$ est
absolument irréductible si cette représentation l'est. \par \vspace{1 mm}
	Soient $F$ un corps, $B \subset F$ un sous-anneau, 
$\rho: G \rightarrow \GL_n(F)$ une représentation de $G$, on suppose de plus
que $\tr(\rho(G)) \subset B$, la discussion ci-dessus implique immédiatement
le:

\begin{cor} \label{ssgenerique}Sous ces hypothèses, si $m$ est un idéal maximal de $B$ de corps
résiduel $k$ de caractéristique $0$, alors la réduction modulo $m$ de
$\tr(\rho)$ est la trace d'une représentation semi-simple $$\rho_m^{ss}: G \rightarrow
\GL_n(\overline{k})$$
$\rho_m^{ss}$ est unique à isomorphisme près, définie sur une extension finie de $k$.
\end{cor}

	Notons qu'il est clair que si $\rho(G) \subset \GL_n(B)$, cela a un
sens de réduire $\rho$ modulo $m$, et qu'alors $\rho_m^{ss}$ est la
semi-simplification de cette réduction. Dans le cas où $B=A(X)$ est une algèbre affinoide, $x \in X$, $m$
l'idéal maximal de $B$ défini par $x$, on notera aussi $\rho_x^{ss}$ pour
$\rho_m^{ss}$. 

\subsubsection{} Supposons $A$ quelconque, mais que pour tout $m \in \Specmax(A)$, la réduction modulo $m$ de
$T$, $T_m: G \rightarrow A/m$ est absolument irréductible. Quitte à
remplacer $A$ par une extension étale finie, il est encore vrai que $T$ est
la trace d'une unique représentation $G \rightarrow \GL_n(A)$ (\cite{Rou}
5.1). \par \vspace{1 mm}
	Supposons finalement que $A$ est intègre, que 
$T: G \rightarrow Frac(A)$ déduit de $T$ soit absolument
irréductible, mais que les $T_m$ ne sont pas tous absolument irréductibles. 
On ne peut plus alors attacher canoniquement de représentation à $T$ mais, au moins sous certaines hypothèses développées plus bas,
un ensemble de représentations non toutes isomorphes en général. Pour nos
applications, le cas intéressant sera celui où $A$ est une algèbre affinoide
intègre de dimension $1$. \par \vspace{2 mm}

\begin{lemme} \label{pseudorig} Soient $F$ un corps local, $X$ un $F$-affinoide
intègre de dimension $1$, et $T: G \rightarrow A(X)$ un pseudo-caractère de
dimension $n$, il existe : \par \vspace{1 mm}
i) Un $F$-affinoide $Y$ régulier intègre, de dimension $1$, fini et surjectif sur
$X$, \par \vspace{1 mm}
ii) Une représentation semi-simple $\rho_{K(Y)}: G \rightarrow \GL_n(K(Y))$
de trace $T$, satisfaisant (ABS). \par \vspace{2 mm}

\end{lemme}

{\it Preuve:} En considérant la composée $T: G \rightarrow A(X) \rightarrow K(X)$,
\S \ref{repattacheesauxpseudorep} assure que $T$ est la trace d'une représentation semi-simple 
$G \rightarrow \GL_n(L)$ satisfaisant (ABS) (cf. \ref{semisimples}), pour une extension
finie $L/K(X)$. Soit $A'$ la normalisation de $A$ dans $L$, 
c'est un anneau de Dedekind car $A$ est intègre noethérien de dimension $1$. 
D'après \cite{BGR} 6.1.2, proposition 4, $A'$ est une $F$-algèbre affinoide
finie sur $A$, on l'écrit $A(Y)$, en particulier $L=K(Y)$. Ceci prouve i) et
ii). 
%La propriété iii) est une conséquence directe de \ref{stable} iv). 
$\square$ \par \vspace{1 mm}

\subsubsection{} \label{euh} Soit $Y$ un $F$-affinoide intègre de dimension 
$1$, $y \in Y(F)$ un point
régulier, $G$ un groupe et $T: G \rightarrow A(Y)$ un pseudo-caractère qui est la trace d'une
représentation semi-simple $G \rightarrow \GL_n(K(Y))$ satisfaisant (ABS). D'après \ref{stable}  iv), $K(Y)^n$ a un $A(Y)_g$ réseau
stable $\Lambda_0$, pour $g\in A(Y)$ ne s'annulant pas en $y$. \par
\vspace{1 mm}
On note $\OOO$ l'anneau local rigide de $A(Y)$ en $y$, 
c'est un anneau de valuation discrète (voir \cite{BGR} 7.3.2 proposition 8), 
on pose $L:=\textrm{Frac}(\OOO)$. 
$\OOO \Lambda_0$ est un $\OOO$-réseau de $L^n:=K(Y) \otimes_{K(Y)} L^n$
stable par $G$. 
On supposera que l'anneau local algébrique en $y$ est principal d'idéal 
maximal engendré par $z \in A(Y)$, on peut toujours
faire cette hypothèse quitte à rétrécir $Y$.\par \vspace{1 mm}

Soient $\Omega_1$ et $\Omega_2$ deux ouverts affinoides de $Y$ contenant $y$,  $\Lambda_i$ un
$A(\Omega)_i$-réseau stable de $L^n$ pour $i=1,2$. 
On dira que $\Lambda_1$ et $\Lambda_2$ 
sont équivalents s'il existe $\Omega_3 \subset \Omega_1\cap \Omega_2$ un 
ouvert affinoide de $Y$ contenant $x$ tel qu'il
existe $r\in \ZZ$ vérifiant $$z^rA(\Omega_3)\Lambda_1=A(\Omega_3)\Lambda_2 \subset L^n$$ On note
$S'$ l'ensemble des classes pour cette relation d'équivalence, $[\Lambda_i]
\in S$ la classe du $A(\Omega_i)$-réseau stable $\Lambda_i$. \par \vspace{1
mm}

\begin{lemme} \label{bij} \begin{itemize}
\item $S$ et $S'$ sont non vides, 
\item L'application $S' \longrightarrow S$, $[\Lambda] \mapsto [\OOO \Lambda]$ est une
bijection,
\item Si $s=[\Lambda] \in S'$, les représentations de $G$ sur les $F$-espaces vectoriels $\Lambda/z\Lambda$ et $(\OOO \Lambda)/z(\OOO \Lambda)$ sont isomorphes. 
\end{itemize}
\par \vspace{1 mm}
Cette dernière classe d'isomorphisme ne dépend pas du choix de $\Lambda$ tel que
$s=[\Lambda]$, on l'appelle la représentation résiduelle de $s$.
\end{lemme}

{\it Preuve:} On a déjà montré le premier point. L'application de l'énoncé est clairement bien définie. Montrons
l'injectivité. Avec les notations du paragraphe ci-dessus, supposons que $\OOO
\Lambda_1=z^p \OOO \Lambda_2$, $p \in \NN$. On peut supposer $p=0$. Si $e^i_1,...,e^i_n$
une $A(\Omega_i)$-base de $\Lambda_i$, on note $M \in \GL_n(\OOO)$ la matrice des
$e^1$ dans la base des $e^2$. On peut trouver $\Omega_3 \subset \Omega_1
\cap \Omega_2$ un affinoide de $Y$ contenant $x$ tel que $M \in
\GL_n(A(\Omega_3))$. Alors $A(\Omega_3)\Lambda_1=A(\Omega_3)\Lambda_2$. \par \vspace{1
mm}
Pour la surjectivité, considérons $\Lambda$ un $\OOO$-réseau stable par $G$ quelconque dans $L^n$. 
On sait que l'image $C$ de $A(Y)[G]$ dans $\End_{K(Y)}(K(Y)^n) \subset
\End_{L}L^n $ est de type fini sur $A(Y)$ par le lemme \ref{stable} i). $C$
est de plus trivialement un sous-$A(Y)$-module de 
$\End_\OOO(\Lambda)$. On choisit
$m_1$,...,$m_r$ une famille $A(Y)$-génératrice de $C$, ainsi que 
$e_1$,...,$e_n$ une $\OOO$-base du $\OOO$-réseau $\Lambda$. On peut trouver un voisinage affinoide $\Omega$ de $x$ dans $Y$ tel que l'ensemble
de tous les $m_i$ soient à coefficients dans $A(\Omega)$ dans la base des
$e_i$. $\oplus_{i=1}^n A(\Omega)e_i$ est alors 
$A(\Omega)$-réseau stable de $L^n$ qui convient. \par \vspace{1 mm}
	Les autres assertions sont immédiates. $\square$ \par \vspace{2 mm}

On rappelle que pour tout affinoide $Y$ réduit, $A(Y)$ est canoniquement
normé par sa norme du sup. Conservant les hypothèses de tout ce paragraphe,
supposant de plus que $G$ est un groupe topologique et que $T: G \rightarrow
A(Y)$ est continue, alors pour tout $A(\Omega)$-réseau stable comme plus
haut, le lemme \ref{stable} v) montre que la représentation déduite $G
\rightarrow \GL_n(A(\Omega))$ est continue. En corollaire du lemme \ref{bij},
on obtient alors le

\begin{cor} \label{continu} Si $T: G \rightarrow A(Y)$ est continue, alors pour tout $\OOO$-réseau stable
$\Lambda$ de $L^n$, la représentation résiduelle de $G$ sur $\Lambda/z\Lambda$ est continue.
\end{cor}

\subsection{Variante d'un lemme de Ribet}

Soit $A$ un anneau de valuation discrète, 
$K$ son corps des fractions, $\p$ un idéal
maximal de $A$, et $k=A/\p$.
Soit $G$ un groupe, $\tau$ un automorphisme de $K$.
Pour $\rho$ une représentation de $G$ sur un espace $V$ de dimension finie
sur $K$,
et $\Lambda$ un $A$-réseau stable (dans la suite nous dirons simplement un 
{\it réseau stable}), on note $\rhob_\Lambda$ la représentation
sur   
$\Lambda/\p\Lambda \simeq k^n$. Par le théorème de Brauer-Nesbitt 
(\cite[30.16]{cr}), $\rhob_\Lambda^{ss}$
ne dépend pas du réseau stable $\Lambda$ (s'il en existe un), on la note
$\rhob^{ss}$.
Pour $\rho$ une représentation sur $K$ ou $k$ on note $\rho^\bot$ la
représentation $g \mapsto \rho(\tau(g))^\ast$.

\begin{prop} \label{ribet2}
Soit $\rho$ une représentation de $G$ de dimension $n$
sur $K$  admettant un réseau stable. On suppose que
$\rho\simeq \rho^\bot$ est absolument irréductible et que 
$\rhob^{ss} \simeq \phi \oplus \phi^\bot \oplus \psi$ où   
$\phi$, $\phi^\bot$ et $\psi$ sont trois
représentations absolument irréductibles deux à deux non isomorphes.

Alors :
\begin{itemize}
\item[a.] Soit il existe un réseau stable $\Lambda$, tel que $\rhob_\Lambda$
admette un sous-quotient $r$ de dimension $2$,
vérifiant $r \simeq r^\bot$ et tel que $r$    
est une extension non triviale de $\phi^\bot$ par $\phi$.
\item[b.] Soit il existe un réseau stable $\Lambda$, tel que
$\rhob_\Lambda \simeq \rhob_\Lambda^\bot$, $\rhob_\lambda$ admet une unique
sous-représentation $r$ de dimension $2$ et un unique sous-quotient $r'$   
de dimension $2$, avec $r$ extension non triviale de $\psi$ par $\phi$,
$r'$
extension non triviale de $\phi^\bot$ par $\psi$, et $r' \simeq r^\bot$.
\end{itemize}
\end{prop}   

\subsubsection{}
Pour prouver la proposition, notons d'abord qu'on peut supposer que $A$
est un anneau de valuation discrète complet, ce que l'on fait.
 En effet, si $A'$ est le 
complété de $A$ en $\p$, $K'$ le corps des fractions de $A'$, $A'$ est de
valuation discrète complet de corps résiduel $k$. 
L'extension $\rho'$ de $\rho$ à $K'$ vérifie encore toutes les hypothèses du
théorème. Si $\Lambda'$ est un $A'$-réseau stable,
$\Lambda = \Lambda'\cap V$ est un  $A$-réseau stable de
$V$, et $\rho_\Lambda \simeq \rho'_{\Lambda'}$, si
bien que la conclusion du théorème pour $K'$ entraîne la conclusion pour $K$. 

\subsubsection{} \label{rappelribellaiche}
Nous prouverons la proposition en utilisant le dictionnaire de
\cite{ribellaiche} et \cite{ribellaicheg}. 
Rappelons quelques notions et résultats utiles de ces articles, auquel nous
renvoyons au lecteur pour de plus amples détails.

Soit $\X$ l'immeuble de Bruhat-tits de $\PGL_n(K)$, $X$ l'ensemble de ses
sommets. Soit $\Sc$ la partie de $\X$ fixe par $\rho(G)$, $S=\Sc \cap X$. 
D'après~\cite[prop. 4.1.1]{ribellaicheg} et \cite[prop. 2.4.1 et 3.3.2]{ribellaiche}, $\Sc$ est contenue dans un sous-appartement $\Ac$
de dimension 2 (car $\rhob^{ss}$ est sans multiplicité, et a trois facteurs
de Jordan-Hölder), clos (\cite[prop.
3.1.3]{ribellaiche}), borné (car $\rho$ est irréductible, \cite[prop
3.2.1]{ribellaiche}), et
les points de $S$ (qui sont donc en nombre fini) correspondent aux classes 
d'homothéties de réseaux stables
par $\rho$ (\cite[prop. 3.1.2]{ribellaiche}). Pour $x \in S$ on note
$\rhob_x$
une réduction $\rhob_\Lambda$ pour $\Lambda$ un représentant de $x$, ce qui
est bien définie à isomorphisme près.

L'ensemble $\Sc$ est un polygone
dans l'appartement $\Ac$ (un plan), dont les sommets sont des sommets de $\Ac$
et les côtés des réunions d'arêtes de $\Ac$ (\cite[2.1]{ribellaicheg}).

À $\Sc$ on peut attacher, comme en \cite[2.1.2]{ribellaicheg} un graphe
orienté d'ensemble de sommets $\{\phi,\phi^\bot,\psi\}$, d'ensemble d'arêtes
$A\subset (\{\phi,\phi^\bot,\psi\}^2-\Delta)$  ($\Delta$ la diagonale)
  connexe en tant que graphe orienté, muni d'une bijection
$c$ de $A$ sur les côtés de $\Sc$, vérifiant les propriétés suivantes
\begin{itemize} 
\item[i)] Pour tout arête $a=(u,v) \in A$ il existe une extension $r_a$
non triviale de $v$ par $u$, et pour $x \in c(a) \cap S$, $r_a$ apparaît
comme sous-quotient de $\rhob_x$. Réciproquement, si $x \in S$, et si $\rhob_x$ contient comme sous-quotient une extension triviale de $v$ par $u$, 
alors $x$ est sur $c(u,v)$. (\cite[prop. 4.2.1]{ribellaicheg})
\item[ii)] 
Si $a$ et $a'$ sont deux arêtes de $A$ non opposées l'une de l'autre,
les côtés $c(a)$ et $c(a')$ du polygone $\Sc$ se coupent en un sommet.
(\cite[prop. 2.2.4]{ribellaicheg})
\end{itemize}

\subsubsection{}
Par hypothèse, il existe un isomorphisme de $K$-espaces vectoriels
$\phi$ de $V$ dans $V^\ast$ tel que 

\begin{eqnarray} \label{entrelac2}
\forall g \in G,\, \rho(\tau(g)) = \varphi^{-1} \rho^\ast(g) \varphi.
\end{eqnarray}

Si $\Lambda$ est un réseau de $V$, on note $\Lambda^\ast$
le réseau dual dans  $V^\ast$.
 L'application $\Lambda \mapsto \Lambda^\ast$ passe
au quotient et définit une bijection naturelle $b$ entre l'ensemble des
sommets de l'immeuble $X$ de $\PGL(V)$ et celui de l'immeuble
$X^\ast$ de $\PGL(V^\ast)$, qui s'étend en un morphisme d'immeubles
de $\X$ dans $\X^\ast$.
Il est clair que si  $\Sc^\ast$ désigne la partie de $\X^\ast$ stable par
$\rho^\ast$, on a $$ S^\ast = b(S).$$
Par ailleurs l'application $\varphi$ induit un
isomorphisme d'immeubles $\varphi_\ast$ de $\X$ dans $\X^\ast$.
On déduit immédiatement
de~(\ref{entrelac2}) que $$\varphi_\ast(S)=S^\ast.$$

L'isomorphisme $t := b^{-1} \varphi_\ast$ de $\X$
laisse donc stables  $\Sc$ et $S$. Si $x \in S$,  
$\varphi$ induit par passage au quotient un isomorphsime
\begin{eqnarray} \label{tbot} \rhob_x^\bot \simeq \rhob_{t(x)}.
\end{eqnarray}

\subsubsection{} Terminons la preuve. Comme le graphe $(V,A)$ est connexe
en tant que graphe orienté, il y a un chemin qui va de $\phi$ à $\phi^\bot$.
Il y a donc deux possibilités :

Soit $(\phi,\phi^\bot)\in A$, dans ce cas
$c(\phi,\phi^\bot)$ est un côté de $\Sc$, et pour 
$x \in c(\phi,\phi^\bot) \cap S$, $\rhob_x$ contient comme sous-quotient
$r_{(\phi,\phi^\bot)}$ qui est une 
extension non triviale de $\phi^\bot$ par $\phi$.
De plus, $t$ envoie côté de $\Sc$ sur côté de $\Sc$, et par~(\ref{tbot}),
les points $x$ du côté $t(c(\phi,\phi^\bot))$ sont tels que $\rho_x$ admette
comme sous quotient
$r_{\phi,\phi^\bot}^\bot$ qui est encore une extension non triviale de 
$\phi^\bot$ par $\phi$. 
On en déduit d'après~\ref{rappelribellaiche} i) que
$t(c(\phi,\phi^\bot))=c(\phi,\phi^\bot)$
et donc, d'après~(\ref{tbot}) que  $r_{(\phi,\phi^\bot)}^\bot\simeq 
r_{(\phi,\phi^\bot)}$, ce qui
montre qu'on est dans le cas a. de la proposition.

Soit $(\phi,\psi)$ et $(\psi,\phi^\bot)$ sont des arêtes de $A$.
Par~\ref{rappelribellaiche} ii), $c(\phi,\psi)$ et
$c(\psi,\phi^\bot)$ se coupent en un sommet $x$. 
Par un raisonnement semblable à celui ci-dessus,
$t(c(\phi,\psi))=c(\psi,\phi^\bot)$ et
$t(c(\psi,\phi^{\bot}))=c(\psi,\phi)$, d'où $t(x)=x$. La réduction $\rhob_x$ vérifie donc $\rhob_x \simeq
\rhob_x^\bot$, et on est dans le cas b. de la proposition.

\section{Déformation $p$-adique de $\chi \oplus 1 \oplus \chi^{\bot}$ }
	
\label{defo}

\subsection{Notations pour les espaces de formes automorphes}
\subsubsection{} On fixe encore $p=v_1v_2$ décomposé dans $E$, 
$N$ un entier premier à $p$, ainsi qu'un isomorphisme de corps $\iota: \CC \rightarrow
\overline{\QQ}_p$, $\CC_p$ le complété de $\overline{\QQ}_p$. \par \vspace{1 mm} La donnée de $v_1$ nous permet d'identifier canoniquement $\U(3)(\QQ_p)$ à
$\GL_3(\QQ_p)$ comme dans \ref{pdeccan}. Soit $K_f$ un compact ouvert de $\U(3)(\A_f)$,  
décomposé place par place, égal à un compact maximal hyperspécial aux places
ne divisant pas $pN\disc(K)$, et à un compact très spécial aux places
divisant $\disc(K)$ mais pas $pN$, et à l'Iwahori $I$ de $\GL_3(\QQ_p)$ 
en $p$. \par \vspace{1 mm} On fixe une représentation
irréductible complexe lisse $J$ de $K_f$ définie sur $\overline{\QQ} \subset
\CC$, triviale restreinte aux places ne
divisant pas $N$. $J$ est de dimension finie, $\iota$ nous permet de la
voir à coefficients dans un corps local fixé $F_0$ et de considérer $J(F)$ pour chaque $F_0 \subset F
\subset \CC_p$ comme étant une $F$-représentation lisse de $K_f$. On note
$\mathcal{B}$ l'espace des fonctions complexes lisses sur
$\U(3)(\QQ)\backslash \U(3)(\A)$, vu comme représentation de $\U(3)(\A)$ par translation à droite. 

\subsubsection{} Si $w=(k_1\geq k_2 \geq k_3) \in \ZZ^3$, l'espace des formes automorphes
pour $\U(3)$ de "poids automorphes" $w$ et de type $(K_f,J)$ est le
$\CC$-espace vectoriel $$S_w(K_f,J,\CC):=\textrm{Hom}_{\U(3)(\RR)\times
K_f}(V_w(\CC) \otimes_{\CC} J,\mathcal{B})$$ 
On note $\HH$ l'anneau engendré par l'algèbre de Hecke globale
(sur $\ZZ$) hors de $pN$ et par l'algèbre d'Atkin-Lehner
$\mathcal{A}(p)$ en $p$. $\HH$ est commutatif, et on notera $\HH_A := \HH
\otimes_{\ZZ} A$, pour tout anneau $A$.
$S_w(K_f,J,\CC)$ est de manière naturelle un module sur $\HH_{\CC}$. 
Si $V$ est un $A[(I\Delta^+I)K_f]$-module, on notera $H^0(V)$ le $\HH_A$-module des
fonctions $$\U(3)(\QQ) \backslash \U(3)(\A_f) \rightarrow V, \, \textrm{telles
que  } \, \, \forall x \in \U(3)(\A_f), \, \, \forall u \in K_f, \, \,
u.f(xu)=f(x)$$
C'est une vérification classique que la donnée de $\iota$ permet de voir 
$H^0(V_{-w}(F_0) \otimes_{F_0} J^*(F_0))$ comme une
$F_0$-structure
de $S_w(K_f,J,\CC)$, et ce comme $\HH_{F_0}$-module. On pose\footnote{Noter le $-w$ dans la définition}
$$S_w^{cl}:=H^0(V_{-w}(F_0)\otimes_{F_0}J^*(F_0))$$ \par \vspace{2
mm}

\subsubsection{} Par la finitude du nombre de classes (\cite{Bor} 5.1), $\U(3)(\A_f)$ s'écrit
comme réunion finie $$\U(3)(\A_f)=\coprod_{i=1}^h \U(3)(\QQ)
x_i K_f$$ de plus $\Gamma_i:=x_i\U(3)(\QQ)x_i^{-1}\cap K_f$ est un groupe fini
car compact discret. Pour un $A[I\Delta^{+}I,K_f]$-module $V$, l'application $A$-linéaire 
$$\varphi_V: H^0(V)
\rightarrow \oplus_{i=1}^h V^{\Gamma_i}, \, \, \, \, f \rightarrow
(f(x_1),...,f(x_h))$$ 
est un isomorphisme
$A$-linéaire fonctoriel en $V$. On déduit que sur les $\QQ$-espaces vectoriels, $V \mapsto H^0(V)$ est exact
et commute à l'extension des scalaires en $A$. \par \vspace{1 mm}
Enfin, on en déduit aussi que $S_w^{cl}$ est de dimension finie sur $F_0$.
De plus, si $V$ est normé par $|.|$ tel que $K_f$ y agisse par automorphismes continus de norme
$1$, $H^0(V)$ est naturellement normé par
$$|f|:=\textrm{Sup}_{x \in \U(3)(\A_f)} |f(x)|=\textrm{Sup}_{i=1}^h|f(x_i)|$$

\subsection{Familles $p$-adiques typées pour $\U(3)$} 

\subsubsection{} \label{notationsfamilles} Soit $w=(k_1 \geq k_2 \geq k_3) \in \ZZ^3$, $\alpha=(\alpha_1,\alpha_2,\alpha_3) \in (\QQ^{\geq 0})^3$, on note
$S_{w}^{cl,\alpha}$ le plus grand sous-$\HH_{\Qpb}$-module de 
$S_w^{cl}\otimes_{F_0}\Qpb$ sur lequel $U_i^w \in \mathcal{A}(p)$
(cf. \ref{versionU}), $1\leq i \leq 3$, n'a que des valeurs propres de valuation
$\alpha_i$. On dira que $w$ est $\alpha$-régulier si $\delta(w) > \alpha_1+\alpha_2-1$.
\par \vspace{1 mm}Si $r \in ]0,1]\cap\QQ$, $x \in \CC_p^3$, on note
$B(x,r)$ la boule fermée de $\CC_p^3$ de centre $x$ de rayon $r$. 

\subsubsection{} Fixons $f \in S_w^{cl} \otimes_{F_0} \Qpb$ une forme propre pour $\HH$. 

\begin{prop} \label{familles} Il existe un corps local $F$, $r \in |F|\cap [0,1[$, $X$ un $F$-affinoide
réduit, $\pi: X \rightarrow B(w,r)$ un $F$-morphisme fini, $a: \HH \rightarrow A(X)^0$ et
$x_f \in X(F)$ tels que : \par \vspace{2 mm}

i) Si $1\leq i \leq 3$, l'application $$\pi^{-1}(\ZZ^3 \cap B(w,r)(F))
\longrightarrow \QQ^{\geq 0}, \, \, \, \, \, x \mapsto
v(a(U_i^{\pi(x)})(x))$$
est constante, disons égale à $\alpha_i$. On pose
$\alpha=(\alpha_1,\alpha_2,\alpha_3)$. \par \vspace{1 mm}
ii) Si $w' \in B(w,r) \cap (w+(p-1)\ZZ^3)$ est $\alpha$-régulier, l'application
$$X(\CC_p) \longrightarrow \Hom(\HH,\CC_p), \, \, x \mapsto (\chi_x: h \mapsto a(h)(x))$$ 
induit une bijection entre $\pi^{-1}(w')$ et l'ensemble des caractères de $\HH$ dans
$S_{w'}^{cl,\alpha}$ comptés sans multiplicités. \par \vspace{1 mm}
iii) $\chi_{x_f}$ est le caractère de $\HH$ sur $f$. \par \vspace{1 mm}
iv) L'image de $\HH$ dans $A(X)^0$ est d'adhérence compacte. \par \vspace{1 mm}

v) La restriction de $\pi$ a chaque composante irréductible de $X$ est surjective
sur $B(w,r)$. Cette propriété est de plus encore vérifiée pour le changement
de base de $\pi$ à tout fermé irréductible de $B(w,r)$. \par \vspace{1 mm}
vi) $\pi^{-1}(B(w,r)(\QQ_p))(\CC_p) \subset X(F)$. \par \vspace{1 mm}
\end{prop}

{\it Preuve:} C'est essentiellement une conséquence des techniques de \cite{Che}, 
à ceci près que nous expliquons ici comment gérer la présence du
type $J$, et le fait que le niveau n'est pas nécessairement net (cf. 4.1.
{\it loc. cit.}). \par
\vspace{1 mm}
Soit $B:=B(0,1)/\QQ_p$ d'algèbre affinoide  $A(B):=\QQ_p<T_1,T_2,T_3>$.
On considère le $A(B)$-module de Banach orthonormalisable
noté $S_{(\tau^{-k_3},\tau^{-k_2},\tau^{-k_1})}(1)$ dans \cite{Che}, \S 3.6.
$M$ est en particulier un $A(B)[I\Delta^{+}I]$-module 
({\it loc. cit.} 3.6.2), $I$ agissant par endomorphismes de norme $1$. Nous ne
prendrons pas exactement la même paramétrisation des poids que dans \cite{Che} (voir
{\it loc. cit.} \S
2.3), nous choisissons ici l'unique telle que les spécialisations $M_{w'}$ en $w' \in
\ZZ^{3,+} \subset B(\QQ_p)$ contiennent $V_{-w'}(\QQ_p)\otimes_{\QQ_p}\tau^{w'-w}$ comme représentation
de $I$ (cf. remarques 3.7.2, {\it loc. cit.} pour la notation $\tau^w$). En
particulier, $M_w$ contient $V_{-w}(\QQ_p)$ comme sous-représentation de $I\Delta^{+}I$. \par \vspace{1
mm}
 On fixe une norme sur $J(F_0)$ et on définit le $A(B)$-module de Banach des formes
$p$-adiques pour $\U(3)$, de type $J$, comme étant $$S:=H^0(M\otimes_{\QQ_p}
J^*(F_0))$$
L'application $\varphi_{S}$ montre que $S \simeq (M \otimes_{\QQ_p} J^*(F_0))^{\Gamma_i}$
est un $A(B)$-isomorphisme. $S$ est donc un $A(B)$-module de Banach
satisfaisant la condition (Pr) de \cite{KB}
\S 2, on a donc bien une théorie des séries caractéristiques pour ses endomorphismes compacts (voir \cite{BMF} A, \cite{KB} \S 2).
\par \vspace{1 mm}
Si $U$ est l'opérateur de Hecke agissant sur $S$ donné par la double classe $I\diag(1,p,p^2)I$, alors $U$ est $A(B)$-linéaire compact (\cite{Che} 3.6.2).
Il faut remarquer de plus que si $w' \in w + (p-1)\ZZ^{3,+}$, $S_{w'} \supset
S^{cl}_{w'}$, et $U_{|S^{cl}_{w'}}$ coïncide alors avec $U_1^wU_2^wU_3^w$. 
On pose $g(T):=\det(1-TU_{|S}) \in 1+TA(B)\{\{T\}\}$, 
c'est une fonction analytique sur $B \times \A^1$. \par \vspace{1 mm} 
On fixe $F \supset F_0$ un corps local contenant les valeurs propres
$\lambda_i$ des $U_i^w$ sur $f$, on note $\alpha_i:=v(\lambda_i)$,
$\alpha_3=0$. Utilisant \cite{Che} 5.4, il est aisé de voir que pour un certain $r \in [0,1[\cap |F|$,
$g$ se restreint en une fonction analytique sur $B(w,r) \times \A^1$ qui se
factorise de la forme $g_{|A(B(w,r))}=PS$, avec $P \in 1+TA(B(w,r))[T]$, $S \in
1+TA(B(w,r))\{\{T\}\}$, $(S,T)=1$ et tels que $\forall x \in B(w,r)(F)$,
$P_x$ a toutes ses racines de valuation $-\alpha_1-\alpha_2$, et $S_x$ n'en a pas.
\par \vspace{2 mm}
	On applique à la donnée ci-dessus la construction des variétés de
Hecke locales faite en \S 6.2 {\it loc. cit.} On découpe tout d'abord un
sous-$\HH_{B(w,r)}$-module $N$ de $S_{B(w,r)}$, localement libre de rang fini sur
$A(B(w,r))$, sur lequel $\det(1-TU)=Q$. Quitte à diminuer $r$, on peut 
remplacer $N$ par un facteur direct $N'$ ayant la propriété que $\forall x \in
B(w,r')(F), \, \, \det(1-T[I\diag(1,1,p)I]_{|N'_x})$ a exactement pour racines
celles de $\det(1-T[I\diag(1,1,p)I]_{|N_x})$ qui sont de valuation $-\alpha_1$.
On définit alors $X$ comme l'affinoide d'algèbre l'image de $\HH_{B(w,r)}$ dans
$\End_{A(B(w,r))}(N')$. Si $X$ n'est pas réduit, on le remplace par sa
nilréduction. On dispose alors de $X \rightarrow B(w,r)$, $a$ et $x_f$ comme dans l'énoncé de la proposition \ref{familles} satisfaisant
i), ii) et iii). Il faut noter que la proposition de 4.7.1 {\it loc. cit.} ("classicité des formes de
poids régulier") est encore vérifiée pour nos espaces de formes typées, la preuve étant
rigoureusement analogue. \par \vspace{1 mm} 
Quitte à remplacer $F$ par une extension finie, vi) vient de ce qu'il n'existe qu'un nombre fini
d'extensions de $\QQ_p$ de degré donné. La première partie de la condition
v) provient de 6.4.2. loc. cit. La seconde se ramène à la première en 
utilisant le lemme 6.2.5 {\it loc. cit.} (bien qu'il soit énoncé avec un idéal $m$ maximal,
une modification immédiate de sa preuve marche pour tout idéal $m$). \par \vspace{1 mm}
	Enfin pour iv), on note que si $0<r'<r \in |F|$, alors $B(w,r') \hookrightarrow B(w,r)$
est compacte sur les algèbres affinoides, ainsi donc que $A(X) \rightarrow A(X)
\widehat{\otimes}_{A(B(w,r))} A(B(w,r'))$. On conclut car $a(\HH) \subset A(X)^0$ et l'image d'un borné par une application compacte
entre espaces de Banach sur un corps local est d'adhérence compacte. $\square$ \par \vspace{2 mm}

%$$\U(3)(\QQ) \backslash \U(3)(\A_f) \rightarrow V_w(F)^* \square J(F)^*$$
%telles que $f(xg)=(g_p^{-1} \times g^{-1})f(x)$.

\subsubsection{} \label{classique} Un point de $X(F)$ sera dit {\it 
classique} si $\pi(x) \in
w+(p-1)\ZZ^{3,+}$ et si le caractère $\chi_x$ de $\HH$ qui lui
correspond par le $ii)$ de la proposition \ref{familles} est réalisé dans $S_{\pi(x)}^{cl}
\otimes_{F_0}F$. Les points classiques sont Zariski-denses dans $X(\CC_p)$
par $v)$, $i)$, $ii)$ ci-dessus ainsi que le lemme 6.4.7 de \cite{Che}. \par
\vspace{2 mm}

\subsection{Pseudo-caractères galoisiens} 

\subsubsection{} \label{repgalx} On se replace dans les hypothèses du paragraphe précédent. 
Fixons $x \in X(F)$ un point classique, 
et $f \in S_{\pi(x)}^{cl}\otimes_{F_0}F$ un vecteur de caractère $\chi_x$
sous $\HH$. On peut considérer un constituant irréductible $\Pi$ de la
représentation automorphe de $\U(3)$ engendrée par $f$, $\Pi_l$ est alors déterminée par
$\chi$ pour tout $l$ premier ne divisant pas $N$. En particulier, la représentation
galoisienne $p$-adique continue semi-simple associée à $\Pi$ dans
\ref{propgal}, que nous noterons disons $V(x)$, ne dépend que de $x$ et 
pas du $\Pi$ choisi, par le théorème de Cebotarev.
Soit $T(x): \Gal \rightarrow F$ la trace de cette représentation, elle est continue et la détermine. \par
\vspace{2 mm}
Si $T: E \rightarrow A(X)$ est une application d'un 
ensemble $E$ à valeur dans $A(X)$, $x \in X(F)$, on notera $T_x$ la
composition de $T$ par l'évaluation en $x$. On renvoie en
\ref{generalitespseudo} pour les généralités sur les pseudo-caractères.

\begin{cor} \label{pseudo} Il existe un unique pseudo-caractère continu 
de dimension $3$, $$T: \Gal_{Np} \rightarrow A(X)$$ 
tel que pour tout $x \in X(F)$ classique, $T_x=T(x)$. Il satisfait $\forall
g \in G$, $T(\tau g \tau)=T(g^{-1})$.
\end{cor}

{\it Preuve:} On définit tout d'abord $T$ sur certaines classes de conjugaisons de
Frobénius. Soit $l=v_1v_2$ un nombre premier totalement décomposé dans $E$
comme en \ref{pdeccan}, on suppose que $(l,Np)=1$, on pose alors
$$T(\Frob_{v_1}):=l^{-1}a([\GL_3(\ZZ_l)\diag(1,1,l)\GL_3(\ZZ_l)])(x)$$ 
$$T(\Frob_{v_2}):=l^{-1}a([\GL_3(\ZZ_l)\diag(l^{-1},1,1)\GL_3(\ZZ_l)])(x)$$ 

Pour voir que c'est bien défini, $X$ étant réduit, il suffit de le vérifier sur les points
classiques $x \in X(F)$, auquel cas c'est juste $T(x)(\Frob_{v_i})$. Par le théorème de 
Cebotarev, la condition $iv)$ de \ref{familles}, et la continuité de $T(x)$
si $x$ est classique, il vient que $T$ se prolonge d'une unique façon sur
$\Gal_{Np}$ en une application continue $K$-valuée. Le fait que c'est un
pseudo-caractère, l'assertion d'unicité, ainsi que $T(\tau g\tau)=T(g^{-1})$,
découlent par Zariski-densité ($X$ étant réduit) en se ramenant aux résultats
analogues pour les $T_x$ avec $x$ classiques qui sont connus.  $\square$ \par
\vspace{2 mm}

\subsubsection{} \label{crisdense} Nous aurons besoin d'un dernier fait:

\begin{prop} Il existe une constante $C>0$ telle que pour tout $x \in X(F)$
classique, si $\delta(\pi(x)) >C$ alors $V(x)$ est cristalline. En
particulier ces points sont Zariski-dense.
\end{prop}

{\it Preuve:} Soit $x \in X(F)$ classique, on choisit $f \in S_{\pi(x)}^{cl} \otimes_{F_0} F$
comme en \S \ref{repgalx}, ainsi que $\Pi$. On sait que $\Pi_p$ est
engendrée par ses $I$-invariants, c'est donc un sous-quotient d'une induite complète du
Borel $\Ind_B(\psi)$ pour un certain $\psi$ comme dans \ref{iwa1}. Pour voir que $V(x)$ est cristalline il
suffit de montrer que $\Pi_p$ est non ramifiée d'après \ref{nonramcris}, 
ou mieux que $\Ind_B(\psi)$ est irréductible. Or on sait que ceci se produit dès que
$\forall i\neq j, \, \, \, \psi_i(p)/\psi_j(p) \neq p$ (\cite{Z} 4.2). Il se trouve que si
$\pi_i(x)=(k_1,k_2,k_3)$, $(a_1,a_2,a_3):=(-k_1-1,-k_2,-k_3+1)$ alors \S
\ref{versionU} montre que
$$p^{a_i}\iota(\psi_i(p))=\frac{\chi_x(U_i^{\pi(x)})}{\chi_x(U_{i-1}^{\pi(x)})}$$
La proposition \ref{familles} i) conclut l'existence de $C$. La seconde assertion s'en déduit à
la manière de \ref{classique} $\square$ \par \vspace{2 mm}

\subsection{Déformations de $\chi \oplus 1 \oplus \chi^{\bot}$}

\subsubsection{}\label{famillechi0} On reprend les notations de
\ref{notationsfamilles}, avec $N:=\cond(\chi_0)$, $K_f = \prod K_l$ où: \par

\begin{itemize}
\item Si est $l$ premier à $p \cond(\chi_0)$, $K_l$ est le sous-groupe défini 
en~\ref{kpdecompose}, 
\ref{kpinerte}, \ref{kpramifiee} selon que $l$ est décomposé, inerte ou
ramifié, \par \vspace{1 mm}
\item Si $l$ divise $\cond(\chi_0)$, $K_l$ est le sous-groupe $K_J(l)$
défini en \ref{types},
$J:=\otimes_{l | N} (J(l) \otimes \chi_{0}\circ \det)$ où $J(l)$ est la représentation de $K_J(l)$ 
définie aussi en \ref{types}.  
\end{itemize}
\vspace{1 mm}
On reconsidère la représentation automorphe $\pi(\chi_0)$ de $\U(3)$,
et on fixe dans tout ce qui suit $\sigma \in \{1,(2,\, 3), (1, \,3, \,2)\}$ 
accessible pour $\pi(\chi_0)$ (cf. \S \ref{accessible}). 
On pose $w:=(\frac{k-1}{2},\frac{k-1}{2},1) \in \ZZ^{3,+}$. 
On peut choisir un $f \neq 0 \in \pi(\chi_0)^I \cap (S_{w}^{cl}\otimes_{F_0}\Qpb)$ 
propre pour $\HH$, de caractère sous $\mathcal{A}(p)$ 
correspondant à $\sigma$ comme dans \ref{accessible}. On applique alors \ref{familles} à ce $f$, puis \ref{pseudo} et
\ref{crisdense} aux conclusions de \ref{familles}, ce qui nous fournit un
corps local $F$, un $F$-affinoide $X$, $x_f \in X(F)$, $B(w,r) \subset \A^3$,
un morphisme $F$-affinoide fini $\pi: X \rightarrow B(w,r)$, un pseudo-caractère $T$,
des $F_i$ et $C>0$ comme dans ces propositions.

\subsubsection{} \label{poids} Il sera commode de raisonner en terme des poids de
Hodge-Tate des $V(x)$ plutôt que de leurs "poids automorphes" $\pi(x)$. On
définit à cet effet (cf. \S \ref{raff-galoisienne})  $\kappa: X \rightarrow
B(\kappa(w),r)$ comme étant la composée de $\pi$ avec $(x,y,z) \rightarrow
(-x-1,-y,-z+1)$. Ainsi, si $x \in X(F)$ est classique, $V(x)$ est de
Hodge-Tate de poids $\kappa(x)$. On pose
$$\kappa_0:=\kappa(w)=(-\frac{k+1}{2},-\frac{k-1}{2},0)$$ \par \vspace{1 mm}

\subsubsection{} Quitte à prendre une extension finie de
$F$, comme le précisera sa preuve, on a la

\begin{prop} \label{bip} Il existe :  \par \vspace{1 mm}

a) Un $F$-affinoide $Y$ intègre régulier de dimension $1$, $y_0 \in Y(F)$,
\par \vspace{1 mm}
b) Une représentation continue semi-simple $$\rho_{K(Y)}: \Gal_{Np} \rightarrow
\GL_3(K(Y))$$  satisfaisant (ABS), $\rho_{K(Y)}^{\bot} \simeq \rho_{K(Y)}$ et
$\tr(\rho_{K(Y)}(\Gal)) \subset A(Y)$, \par \vspace{1 mm}
c) Un $F$-morphisme $$\kappa=(\kappa_1,\kappa_2,\kappa_3): Y
\longrightarrow \A^3, \, \, \kappa_3:=0, \, \, \kappa(y_0)=\kappa_0,$$ \par
d) Une partie $Z \subset Y(F)$ telle que $\kappa(Z) \subset
\kappa_0+(p-1)\ZZ^{3,--}$, \par \vspace{1 mm}

e) Des fonctions $F_1$, $F_2$ et $F_3$ dans $A(Y)$, chacune de valuation constante
sur $Y(\CC_p)$, \par \vspace{2 mm}

\noindent le tout satisfaisant aux propriétés suivantes : \par \vspace{2 mm}

i) Pour tout ouvert affinoide $\Omega$ de $Y$ contenant $y$, 
la fonction $$x \in Z \mapsto -\delta(\kappa(x)) \in \NN$$ 
est non majoré sur $Z \cap \Omega$, d'image incluse dans $\NN^{\geq \textrm{Max}(C,k)}$. En particulier,
$Z \cap \Omega$ est Zariski-dense dans $\Omega$, \par \vspace{1 mm}

ii) Si $z \in Z \cup \{y_0\}$, $\rho_z^{ss}:=\rho_{K(Y),z}^{ss}$ (cf.
\ref{ssgenerique}) est la représentation galoisienne attachée à une représentation automorphe $\Pi$ irréductible de $\U(3)$ telle que
$Hom_{K_f}(J,\Pi) \neq 0$, \par \vspace{1 mm}
iii) $\rho_{y_0}^{ss} \simeq 1 \oplus \chi_p \oplus
\chi_p^{\bot}$. \par \vspace{1 mm}
iv) Si $z \in Z \cup \{y_0\}$, $(\rho_z^{ss})_{|D_{v_1}}$ est cristalline
de poids de Hodge-Tate $\kappa(z)$. Elle est raffinée par
$$(p^{\kappa_1(z)}F_1(z),p^{\kappa_2(z)}F_2(z),p^{\kappa_3(z)}F_3(z))$$ \par
\vspace{1 mm}
v) Ce raffinement est $\mathcal{R}(\sigma)$ en $y_0$.
\end{prop}

{\it Preuve:} Soit $B \subset B(\kappa_0,r) \subset \A^3$, le fermé affinoide de
$B(\kappa_0,r)$ défini par $x_3=0$, et $x_2=2 x_1+\frac{k-3}{2}$. $\kappa_0
\in B(F)$. On pose:

$$\mathcal{Z}:=\{ z \in B(F) \cap (\kappa_0+(p-1)\ZZ^{3,--}), \, \,
-\delta(z) > \textrm{Max}(C,k,\alpha_1+\alpha_2-1) \}$$

Le choix de $B$, assez arbitraire, est tel que les fonctions $x_2-x_1$ et $-x_2$ sont non bornées sur 
$U \cap \mathcal{Z}$ pour tout $U$ ouvert affinoide de $B$ contenant $\kappa_0$.
\par 
On considère $X_B:= X \times_{B(\kappa_0,r)} B$, c'est un $F$-affinoide de dimension
$1$, $x_f \in X_B(F)$. Le morphisme déduit par extension des scalaires $\kappa_B: X_B \rightarrow B$
est encore fini, surjectif restreint à chaque composante irréductible de
$X_B$ d'après \ref{familles} v). On choisit alors 
$X'$ une composante irréductible (réduite) de dimension $1$ de $X_B$
contenant $x_f$. \par 
	Le pseudo-caractère $T$ peut être vu à valeur dans $A(X')$, par
composition $A(X) \rightarrow A(X_B) \rightarrow A(X')$, et on peut
appliquer le lemme \ref{pseudorig} à la donnée de $A(X')$ 
et $T$. Il nous fournit un $F$-affinoide intègre $Y$, régulier de dimension $1$, muni
d'un morphisme fini et surjectif $h: Y \rightarrow X'$, ainsi qu'une représentation semi-simple 
$\rho: G \rightarrow \GL_3(K(Y))$ de trace $T$, satisfaisant (ABS). Quitte à remplacer $F$ par une extension finie, 
on peut choisir $y \in Y(F)$ tel que $h(y)=x$. Notant que
$T(g^{-1})=T(\tau.g.\tau)$ d'après \ref{pseudo}, on a prouvé a) et b). \par 

On définit $\kappa$ comme étant la composée $Y \overset{h}{\rightarrow} X'
\hookrightarrow X_B \overset{\kappa_B}{\rightarrow} B$. Notons que $\kappa: Y \rightarrow B$ est plat, car c'est le cas des
extensions finies d'anneaux de Dedekind. Cela prouve c). On pose $Z:=\kappa^{-1}(\mathcal{Z})$, il satisfait d) par
\ref{poids} et i) par choix de $\mathcal{Z}$ et platitude de $\kappa$. 
$X'$ étant un fermé de $X$, on peut y restreindre les $F_i$ de \ref{famillechi0}, et les définir sur $Y$
en les composant au morphisme $h: Y \rightarrow X'$. Ce sont ces derniers que l'on
choisit pour e), l'assertion sur les valuations des $F_i$ est déjà satisfaite
sur $X$ par construction (cf. \ref{familles}). \par 
	Soit $z \in Z \cup \{y_0\}$, $\rho_{K(Y),z}^{ss}$ est la représentation
semi-simple de $\Gal_{Np}$ de trace $T_z=T_{h(z)}$ (cf. \ref{ssgenerique},
\ref{pseudo}). Par choix de $Z$ et \ref{familles} iii), $h(z) \in X(F)$ 
est un point classique, et l'assertion $ii)$ découle donc de
\ref{repgalx} et \ref{pseudo}. On déduit alors iii) de ii), \ref{familles}
iii) et du choix de $f$ dans \ref{famillechi0}. \par 
	La première assertion de iv) est alors une conséquence de ii),
\ref{raff-galoisienne}, \ref{poids}, ainsi que \ref{crisdense} et le fait
que $-\delta(\kappa(z))>C$ par i). La seconde assertion, ainsi que v),
découlent de \ref{accessible} et \ref{raff2}. $\square$ \par \vspace{2 mm}

{\it Remarques:} On rappelle que toute la construction ci-dessus dépend du choix
initial du $\sigma \in \got{S}_3$ accessible pour $\pi(\chi_0)$. Ceci fait, le second choix réellement 
effectué dans la construction ci-dessus est celui de la composante irréductible $X'$ de
$X'_B$ passant par $x_f$. Il semble difficile d'évaluer le nombre de
composantes irréductibles de $X$ (ou de $X_B$) au voisinage de $x$. En ce qui concerne ce texte,
chaque choix de composante permet de conclure dans la section \ref{fin}. \par \vspace{1 mm}

\section{Construction de l'extension} \label{fin}

On reprend les notations de la proposition \ref{bip}, où l'on a fixé $\sigma=(2,
3)$.

\subsection{Irréductibilité générique}

\begin{prop} $\rho_{K(Y)}$ est absolument irréductible. \label{irrgen}
\end{prop}

{\it Preuve:} D'après \ref{bip} b), $\rho_{K(Y)}$ vérifie la propriété (ABS) 
(voir ~\ref{semisimples}). Il suffit donc de montrer que $\rho_{K(Y)}$ est 
irréductible. Supposons par l'absurde que $\rho_{K(Y)} \simeq \rho_{1,K(Y)} \oplus
\rho_{2,K(Y)}$, $\rho_{i,K(Y)}$ étant une $K(Y)$-représentation de $G$ de dimension $\neq
0$. Le lemme \ref{stable} ii) montre que $\tr(\pi_i(G)) \subset A(Y)$, car
$A(Y)$ est régulier d'après \ref{bip} a). Soit $z \in Z$,
l'évaluation en $z$ de $\tr(\rho(g)))=\tr(\rho_1(g))+\tr(\rho_2(g)), g \in
\Gal_{Np}$ a un sens et montre que $\rho_{z}^{ss}$ est réductible (cf.
\ref{ssgenerique}). Nous allons montrer que c'est absurde par notre choix du raffinement. \par 
	D'après \ref{bip} e), $v(F_i(.)): X(F) \rightarrow \QQ, x \mapsto
v(F_i(x))$ est constante, on la note $\alpha_i$. En évaluant en $y_0$, 
\ref{bip} iv), v) ainsi que \ref{exemple} appliqué à $(3,\ 2)$ montrent que:
$$(\alpha_1,\alpha_2,\alpha_3)=(1,\frac{k-1}{2}, -\frac{k+1}{2})$$
Notons qu'avec ce choix, $\forall i,j \in \{1,2,3\}, \, \alpha_i \neq 0
\textrm{  et  } \, \alpha_i+\alpha_j \neq 0$. De plus, si $i \neq j$, $|\alpha_i| < k$ et $|\alpha_i+\alpha_j|<k$. \par
	D'après \ref{bip} iv), $\rho_z^{ss}$ est cristalline de poids de Hodge-Tate $k_i:=\kappa_i(z)$,
avec $k_1 < k_2 < k_3$, et les valeurs propres de son Frobénius cristallin
ont pour valuation $k_1+\alpha_1$, $k_2+\alpha_2$ et $k_3+\alpha_3$. 
Par le théorème $A$ de \cite{CF}, pour voir que $\rho_z^{ss}$ est irréductible, il suffit de voir
que $D:=D_{cris}(\rho_z^{ss})$ n'admet pas de sous-module filtré faiblement
admissible. Si $D'$ est un tel sous-module de rang $1$, alors par faible
admissibilité $t_H(D')=t_N(D')$  (cf. {\it loc. cit.}) entraîne $k_i=k_j+\alpha_j$ pour un couple
$(i,j)$. Mais par \ref{bip} i), $|k_i-k_j| > k$, les inégalités sur les
$\alpha_i$ entraînent donc $i=j$ puis $\alpha_i=0$, ce qui est absurde. 
De même, si $D'$ est un sous-module filtré faiblement admissible de rang $2$ de $D$, 
on arrive à une contradiction en résolvant $k_{i'}+k_{j'}=k_i+k_j+\alpha_i+\alpha_j$. $\square$
\par \vspace{2 mm}

{\it Remarques:} Une légère modification de la preuve montrerait qu'en fait
${\rho_{K(Y)}}_{|D_{v_1}}$ est irréductible. Notons de plus que l'argument
d'irréductibilité étant local en $p$, il n'utilise pas le fait que les
représentations galoisiennes attachées aux représentations automorphes 
stables tempérées sont globalement irréductibles, mais
simplement~\ref{nonramcris}. Le point clé est que nous disposons d'un $\sigma$
accessible tel que le raffinement $\mathcal{R}(\sigma)$ est aussi éloigné
que possible du raffinement ordinaire au point $y_0$.

\subsection{L'inertie aux places ne divisant pas $p$}
\subsubsection{}
Nous allons préciser l'action de
l'inertie aux places de $E$ ne divisant pas $p$. Pour énoncer
commodément nos résultats, on introduit 
$$\rho'_{K(Y)}:=\rho_{K(Y)}\otimes_F (\chi_p^{\bot})^{-1}.$$

\begin{prop} \label{actin}
Soit $w$ une place de $E$ au-dessus de $l \not = p$.
Alors \begin{itemize}
\item Si $w$ ne divise pas $\cond(\chi_0)$, $\rho_{K(Y)}$ et $\rho'_{K(Y)}$
sont non ramifiées en $w$.
\item Si $w$ divise $\cond(\chi_0)$, on a
$$\dim_{K(Y)}(\rho'_{K(Y)})^{I_w} = 2.$$ 
\end{itemize}
\end{prop}   

Le reste de la sous-section est consacré a la preuve de cette proposition.
D'après \ref{bip} a), b) et \ref{stable} iv) appliqué en l'idéal maximal de
$y_0$, on peut trouver $g \in A(Y)$ avec $g(y_0) \neq 0$, tel que $\rho_{K(Y)}$
admette un $A(Y)_g$-réseau stable. On note  $\rho$ la représentation de $G$ définie par ce réseau, 
$S$ l'ensemble fini des zéros de $g$, et pour $y \in Y
\backslash S$, $\rho_y$ la réduction de $\rho$ en $y$. \par
Remarquons que $\rho_y$ est bien définie à isomorphisme près, et non plus seulement à semi-simplification près.  
De plus, $\rho_{K(Y)}$ étant semi-simple, $\rho_z$ l'est aussi
pour un sous-ensemble infini de $Z \cap (Y \backslash S)$, que l'on note $Z'$. 
Enfin, si $z \in Z'$, alors $\rho_z^{ss} \simeq \rho_z$
est par \ref{bip} ii) la représentation galoisienne attachée à une
représentation automorphe notée $\Pi(z)$ de $\U(3)$, comme en \ref{repgalx}.

\subsubsection{Le cas $w$ ne divisant pas $\cond(\chi_0)$}

Il suffit de montrer que les $\rho_z$ sont non ramifiées pour $z \in Z'$.
 Par construction (cf. \ref{famillechi0}), $\Pi(z)$ a un vecteur fixe par le 
compact maximal $K_l$. Si $l$ est non ramifié dans $E$, $K_l$ est 
hyperspécial (\ref{kpdecompose}, \ref{kpinerte}) et $\Pi(z)$ 
est non ramifié, si bien que la représentation galoisienne $\rho_z$ associée
l'est aussi, d'apres la propriété ~\ref{etoile} \S \ref{propgal}. 
Si $l$ est ramifié, $K_l$ est très spécial (\ref{kpramifiee}),
si bien que d'après \cite{labesse} p.88, le changement de base $\pi_E$
est non ramifié, et on en déduit encore que $\rho_z$ est non ramifiée,
cette fois d'après la proposition~\ref{langram}.

\subsubsection{Le cas $w$ divisant $\cond(\chi_0)$}

Par construction
(cf. \ref{famillechi0}), on a  $$\textrm{Hom}_{K_J(l)}(J(l),(\Pi(z) \otimes
(\chi_{0}
\circ \det)^{-1})_l) \neq 0$$ D'après
\ref{types}, il existe un sous-groupe ouvert $I'_w$ de $I_w$ tel que pour
tout $z \in Z'$ on a $\rho'_{z}(I'_w)=1$. On en tire $\rho'_{K(Y)}(I'_w)=1$.
\par
Notons que pour démontrer la proposition, il suffit de le faire après une extension
finie de $K(Y)$. Mais il existe une extension finie $F'/F$ telle que 
la représentation ${\rho'_{K(Y)}}_{|I_w}$, qui se factorise par le groupe {\it fini}
$I_w/I'_w$, soit définie sur $F'$. Autrement dit, si $L$ est une extension
composée de $K(Y)$ et $F'$ et $\rho'_L := \rho'_{K(Y)}\otimes_{K(Y)} L$, 
alors ${\rho'_{L}}_{|I_w}$ est isomorphe à $\theta \otimes_{F'} L$, où $\theta$ est une
représentation de $I_w$ sur $F'$ triviale sur $I'_w$. 
Comme $\theta$ est nécessairement semi-simple, on en déduit par évaluation
des traces en $y_0$ que $(\rho'_{|I_w})^{ss}_{y_0} \simeq \theta$. En particulier, $\theta \simeq 1 \oplus 1 \oplus
((\chi_p^\bot)^{-1})_{|I_w}$ d'après~\ref{bip} point iii). La proposition en découle.

\subsection{Application des méthodes à la Ribet et Kisin}

Dans ce paragraphe on tire les fruits de la variante du lemme de Ribet 
démontré dans la proposition \ref{ribet2}, et la proposition \ref{kisin} "à la Kisin".
Introduisons encore quelques notations : on pose 
$$u:=\chi^\bot(\Frob_{v_1})$$ Notons alors qu'on a $\chi(\Frob_{v_1})=up^{-1}$.
Notons $\DD$ le foncteur $\DD(V):=D_{cris}(V_{|D_{v_1}})$, qui est exact à 
gauche. L'action du Frobenius cristallin
$\varphi$ sur les droites 
$\DD(\chi_p^\bot)$, $\DD(\chi_p)$ et $\DD(1)$ est la multiplication 
respectivement par $u$, $up^{-1}$ et $1$ (cf. \S \ref{notationshodge}). Ces trois nombres sont deux à 
deux distincts, puisque leur valuations (respectivement $-(k-1)/2$, $-(k+1)/2$
et $0$) le sont.

\begin{prop} \label{bientotlafin}
Il existe une représentation continue $\rhobar: \Gal \rightarrow \GL_3(F)$
vérifiant
\begin{itemize}
\item[i.] Pour toute place $w$ de $E$ ne divisant pas $p$ on a
\begin{eqnarray}
\dim_{F}(\rhobar  \otimes (\chi_p^\bot)^{-1})^{I_w} &\geq& 2 \text{ si }
w | \disc(\chi_0) \\ \dim_{F}(\rhobar  \otimes (\chi_p^\bot)^{-1})^{I_w}&= & 3 \text{ si } w \not| \disc(\chi_0)
\end{eqnarray}
\item[ii.] $\DD(\rhobar)^{\phi=u}$ est non nul.
\item[iii.] On a $\rhobar^{ss} \simeq \chi_p \oplus \chi_p^\bot \oplus 1$.
Une des deux assertions suivantes est vraie :
\begin{itemize}
\item[a.] Soit $\rhobar$
admet un sous-quotient $r$ de dimension $2$,
vérifiant $r \simeq r^\bot$ et tel que $r$    
est une extension non triviale de $\chi_p^\bot$ par $\chi_p$.
\item[b.] Soit 
$\rhob \simeq \rhob^\bot$, $\rhob$ admet une unique
sous-représentation $r_1$ de dimension $2$ et un unique sous-quotient $r_2$   
de dimension $2$, avec $r_1$ extension non triviale de $1$ par $\chi_p$,
$r_2$ extension non triviale de $\chi_p^\bot$ par $1$, et $r_1 \simeq {r_2}^\bot$.
\end{itemize}
\end{itemize}
\end{prop}
{\it Preuve :}
Notons $\anneau$ l'anneau local rigide de $Y$ en $y_0$, $L$ le
corps des fractions de cet anneau, et 
$\rho_L := \rho_{K(Y)}\otimes_{K(Y)} L$. L'anneau $\anneau$ est
de valuation discrète, de corps résiduel $F$. Le représentation $\rho_L$ est irréductible d'après la proposition~\ref{irrgen}

D'après la proposition~\ref{bip} ii), $\overline{\rho_l}^{ss}$ est la somme de trois caractères, $\chi_p$, $\chi_p^\bot$ et $1$. Ces trois caractères sont 
deux à deux distincts  (ils n'ont pas les mêmes poids) et on est donc en mesure d'appliquer la proposition~\ref{ribet2} à $\rho_L$.
Cette proposition affirme précisément 
l'existence d'un $\anneau$-réseau $\Lambda \subset L^3$ stable par 
$\rho_L$, tel que la représentation réduite associée 
$\rhobar:=\overline{\rho_L}_\Lambda$ vérifie soit la condition iii.a soit la 
condition iii.b de~\ref{bientotlafin}. $\anneau$ étant de valuation
discrète, il résulte immédiatement de la proposition~\ref{actin} que $\rhobar$ vérifie la propriété ii.

Nous allons montrer que $\rhobar$ vérifie ii. Le lemme~\ref{bij} appliqué
à $\rho_L$ a $y_0$ et à la classe d'homothéties $s$ du réseau $\Lambda$ 
donne l'existence d'un ouvert affinoide $\Omega \subset Y$ 
contenant $y_0$, telle que la représentation $\rho_{L}$ admet un 
$A(\Omega)$-réseau stable $M$ ; notant $\rho$ la représentation $\Gal \rightarrow \Gl(M)$,  le lemme~\ref{bij} 
assure que $\rho_{y_0} \simeq \rhobar$. $\rhobar$ est continue d'après
le corollaire \ref{continu}.

Nous allons maintenant appliquer la proposition~\ref{kisin} à
la donnée de $\rho_{|D_{v_1}}: D_{v_1} \rightarrow \GL_3(A(\Omega))$, de
$\kappa_{|\Omega}$ et des ${F_i}_{|\Omega}$ (cf. \S
\ref{generalitesfamilleslocales}). On choisit pour ensemble "$Z$" {\it loc. cit.}
l'ensemble $Z \cap \Omega$ auquel on enlève le sous-ensemble des $z$
tels que $\rho_z$ n'est pas semi-simple. Ce dernier est fini car
$\rho_{K(\Omega)}$ est semi-simple. La proposition \ref{bip} assure que
les hypothèses i) à vi) de \ref{generalitesfamilleslocales}
de la proposition~\ref{kisin} sont vérifiées, iii) par notre
choix de $Z$. On voit donc que 
$$D_{cris,v_1}(\rho_{y_0})^{\phi=F_1(y_0)\kappa_1(y_0)}\text{ est non nul},$$
 où  $\kappa_1$
et $F_1$ sont ceux donnés par la proposition~\ref{bip}. D'après les points iv)
et v) de \ref{bip}, $p^{\kappa_1(y_0)}F_1(y_0)$ est la première 
valeur propre
du raffinement $\mathcal{R}((3,2))$ de ${\rho_{y_0}}_{|D_{v_1}}$, donc 
d'après~\ref{exemple}, $p^{\kappa_1(y_0)}F(y_0)=\chi^\bot(p)=u$, ce qui 
prouve ii. $\square$.

\subsection{Élimination du cas a.}

Nous voulons montrer, par l'absurde, que l'on n'est pas dans le cas iii. a. de la proposition précédente. 
On se place donc dans ce cas. La représentation $\rhobar$ admet comme 
sous-quotient une extension non triviale $r$ de $\chi_p^\bot$ par $\chi_p$. Par conséquent, 
$\rhobar':=\rhobar\otimes (\chi_p^{\bot})^{-1}$ contient comme sous-quotient
$r':=r \otimes (\chi_p^\bot)^{-1}$, extension non triviale de $F$ 
(la représentation triviale sur $F$) par $F(1)$ (le caractère cyclotomique sur $F$).

\begin{lemme} \label{lemmecris}
La représentation $r'$ est cristalline en $v_1$ et en $v_2$.
\end{lemme} 
{\it Preuve :} 
Il suffit de le prouver pour $r$, car $\chi_p^\bot$ est cristallin en $v_1$ et $v_2$, $\chi^\bot$ étant non ramifié en ces places. 
De plus, comme $r\simeq r^\bot$, il suffit de le prouver en $v_1$.
Comme $\DD$ est exact à gauche,  ainsi que le foncteur 
$V \mapsto \DD(V)^{\phi=u}$, on voit que 
$$\dim_F \DD(\rhobar)^{\phi=u} \leq \dim_F 
\DD(r)^{\phi=u} + \dim_F 
\DD(1)^{\phi=u}.$$
Comme $\DD(1)^{\phi=u}=0$ car $u \not=1$, il résulte de i. que 
$\DD(r)^{\phi=u}$ est non nul. 
 
Utilisant encore que $\DD$ est exact à gauche il vient 
$$\DD(\chi_p) \subset \DD(r)$$ et il y a donc dans $\DD(r)$ deux 
droites sur lesquelles $\varphi$ agit par $u$ et par $up^{-1}$ ce qui implique
qu'elles sont distinctes. On en déduit que $\dim_F \DD(r)=2$, i.e. $r$ est cristalline en $v_1$.
$\square$

\begin{lemme} La représentation $r'$ est non ramifiée en toutes les places $w$ ne divisant pas $p$.
\end{lemme}
{\it Preuve :} Si $w$ ne divise pas $\cond(\chi_0)$, $\rhobar'$ est non 
ramifiée en $w$ d'après ii., et $r'$ non plus.

Si $w$ divise $\cond(\chi_0)$, l'exactitude à gauche du foncteur des 
invariants sous $I_w$ donne $$\dim_F (\rhobar')^{I_w} \leq {\dim_F r'}^{I_w} + 
\dim_F ((\chi_p^\bot)^{-1})^{I_w}$$ Comme $((\chi_p^\bot)^{-1})^{I_w}=0$, il découle de
i. que $r'$ est non ramifié. $\square$. \par \vspace{2 mm}

L'existence de $r'$ est alors en contradiction avec le lemme
bien connu suivant, 

\begin{lemme} Soit $E$ est un corps quadratique imaginaire, $F/\QQ_p$ un corps
local, alors il n'existe pas d'extension non triviale de représentations continues de $\Gal$ 
de $F$ par $F(1)$ qui soit non ramifiée hors de $p$ et cristalline en les places divisant $p$.
\end{lemme}

{\it Preuve:} 
La théorie de Kummer entraine que la dimension sur $F$ du groupe de ces extensions 
est égale au rang de $\OOO_{E}^*$ 
(voir par exemple \cite{rubinlivre}, proposition
1.6.4). Ce dernier groupe est fini si $E$ est quadratique imaginaire. $\square$ \par
\vspace{2 mm}

\begin{remarque}  \label{corps CM}
Dans le lemme précédent, il est essentiel que les extensions considérées 
soient non ramifiées (resp. cristallines) à {\it toutes} les 
places. Si on relâche la condition à une place quelconque, de telles
extensions non triviales existent. Par ailleurs il est essentiel que le corps
de base $E$ soit quadratique imaginaire. Pour tout autre corps, à part $\Q$,
l'énoncé précédent serait mis en défaut : on construirait par la théorie de 
Kummer une extension non triviale, en partant d'une unité de $E$ qui n'est 
pas racine de $1$. En fait, l'énoncé précédent correspond, par les 
conjectures de Bloch-Kato, au fait que la fonction $\zeta$ du corps $E$ ne 
s'annule pas en $s=0$. Comme on le sait, ceci n'est vrai que pour $E$ 
quadratique imaginaire ou $\Q$.
C'est ici, et ici seulement à part la formule de multiplicité~\ref{endo-1}, 
que l'on utilise dans ce papier l'hypothèse que $E$ est quadratique 
imaginaire. Tout le reste marcherait tout aussi bien
en travaillant avec $E/F$ extension $CM$, et le groupe unitaire $\U(3)$ 
compact à toutes les places à l'infini sur $F$. Dans ce cadre, la preuve
de~\ref{endo-1} montre que l'analogue de la représentation automorphe 
$\pi(\chi_0)$ existe si et seulement si 
$\varepsilon=(-1)^{\dim_\Q F}$. En particulier, si $\dim F$ est pair, on peut avoir $L(\chi_0,1/2) \not=0$, et on ne
s'attend alors 
pas à ce qu'une extension de $\chi_p$ par $1$ cristalline existe. 
On voit que dans ce cas, c'est un élément non trivial de $H^1_f(E,\QQ_p(1))$ 
qu'aurait construit notre méthode. 
\end{remarque}
\subsection{Fin de la preuve}

On est donc dans le cas b. de la proposition. 
\begin{lemme} $r_1$ et $r_2$ sont cristallines en $v_1$ et en $v_2$.
\end{lemme}
{\it Preuve :} Comme $r_1 \simeq r_2^\bot$, la cristallinité de $r_1$ en
$v_2$  (resp. $v_1$) équivaut à celle de $r_2$ en $v_1$ (resp. $v_2$). Il suffit donc de prouver que $r_1$ et
$r_2$ sont cristallines en $v_1$. 

En ce qui concerne $r_1$, extension non triviale de $1$ par $\chi_p$, on 
observe que les poids de Hodge-Tate en $v_1$
de $\chi_p$ et $1$ sont respectivement $-(k+1)/2$ et $0$, avec 
$-(k+1)/2 \leq -2 =0-2 $. Une telle extension est automatiquement cristalline
d'après \cite{FP} proposition 3.1.

Considérons la représentation $r_2$, pour laquelle le même argument de poids 
ne marche plus. Cependant la cristallinité de $r_2$ en $v_1$ se montre exactement
 comme celle de $r$ dans le lemme~\ref{lemmecris}, à l'aide du point i.
$\square$.

La représentation $r_1$ fournit une extension de $1$ par $\chi_p$, non 
triviale, qui a bonne réduction aux deux places divisant $p$ d'après le 
lemme précédent. Comme $\chi_p$ n'est pas le caractère cyclotomique, une
telle extension a automatiquement bonne réduction (voir l'introduction) aux places ne divisant 
pas $p$ : cf. par exemple~\cite[lemme 1.3.5]{rubinlivre}. Ceci prouve le 
théorème~\ref{principal}.

\bigskip
\bigskip
\tiny
\noindent
Joël Bellaïche : Institute for Advanced Studies \\ Einstein drive, 
Princeton, NJ 08540, États-unis \\ email : jbellaic@ias.edu
 \\ {  }\\ 
Gaëtan Chenevier : Département de mathématiques et applications, École normale supérieure \\45 rue
d'Ulm, 75005 Paris, France
 \\email : chenevie@ens.fr

\end{document}